\documentclass[12pt]{article}
\usepackage{amscd,amsfonts,amssymb,amsmath,latexsym,array,hhline}

\usepackage{theorem}

\renewcommand{\bar}{}
\renewcommand{\tilde}{}

\renewcommand{\rho}{r}

\newcommand{\LL}{{\mathbb L}}

\newcommand{\NN}{{\mathbb N}}
\newcommand{\RR}{{\mathbb R}}

\newcommand{\eps}{{\varepsilon}}

\renewcommand{\phi}{\varphi}

\numberwithin{equation}{section}

\renewcommand{\section}[1]{\refstepcounter{section}
                            {\noindent\large\bf\thesection. #1}}

\makeatletter

\renewcommand{\section}{\@startsection{section}{1}{0pt}{30pt}{6pt}{\large\bf}}
\def\dot{\hspace{-16pt}.\hspace{-2pt} }

\renewcommand{\@makefnmark}{}
\renewcommand{\@cite}[2]{[{#1\if@tempswa ; #2\fi}]}

{\hspace*{.1pt}\hspace*{\fill}$\square$\vskip\baselineskip}
\newtheorem{theorem}{Theorem}[section]
\newtheorem{lemma}[theorem]{Lemma}

\usepackage{color}

\makeatother

\begin{document}

\title{Perpetual American options in diffusion-type \\ models with running maxima and drawdowns}

\author{Pavel V. Gapeev\footnote{London School of Economics,
Department of Mathematics, Houghton Street, London WC2A 2AE, United
Kingdom; e-mail: p.v.gapeev{\char'100}lse.ac.uk}
\and Neofytos Rodosthenous\footnote{Queen Mary University of London, School of Mathematical Sciences, Mile End Road, London E1 4NS,
United Kingdom; e-mail: n.rodosthenous{\char'100}qmul.ac.uk}}
\maketitle


\begin{abstract}
We study perpetual American option pricing problems in an extension of the Black-Merton-Scholes model
in which the dividend and volatility rates of the underlying risky asset depend on the running values
of its maximum and maximum drawdown. The optimal exercise times are shown to be the first times
at which the underlying asset hits certain boundaries depending on the running values
of the associated maximum and maximum drawdown processes. We obtain closed-form solutions
to the equivalent free-boundary problems for the value functions with smooth fit at the
optimal stopping boundaries and normal reflection at the edges of the state space
of the resulting three-dimensional Markov process.
The optimal exercise boundaries for the perpetual American options on the maximum of the market depth with fixed and
floating strikes are determined as the minimal solutions of certain first-order nonlinear ordinary differential equations.
\end{abstract}




\section{\dot Introduction}

       The main aim of this paper is to present closed-form solutions
       to the discounted optimal stopping problem of (\ref{V5b}) for the
       running maximum $S$ and the running maximum drawdown $Y$ associated with
       the continuous process $X$ defined in (\ref{dX4})-(\ref{S4}).
              This problem is related to the
       option pricing theory in mathematical finance, where
       the process $X$ can describe the price
       of a risky asset (e.g. a stock) on a financial market.
       The value of (\ref{V5b}) can therefore be interpreted as the
       rational (or no-arbitrage) price of a perpetual American option in a
       diffusion-type extension of the Black-Merton-Scholes model
       (see, e.g. Shiryaev \cite[Chapter~VIII; Section~2a]{FM}, Peskir and Shiryaev
       \cite[Chapter~VII; Section~25]{PSbook}, and Detemple \cite{Det},
       for an extensive overview of other related results in the area).

       Optimal stopping problems 
       for running maxima of some diffusion processes given linear costs were studied by Jacka \cite{Jmax},
       Dubins, Shepp, and Shiryaev \cite{DSS}, and Graversen and Peskir \cite{GP1}-\cite{GP2}
       among others, with the aim of determining the best constants in the corresponding maximal inequalities.
       A complete solution of a general version of the same problem was obtained in Peskir \cite{Pmax},
       by means of the established maximality principle, which is equivalent to the superharmonic
       characterisation of the value function. Discounted optimal stopping problems for certain
       payoff functions depending on the running maxima of geometric Brownian motions were
       initiated by Shepp and Shiryaev \cite{SS1}-\cite{SS2} and then considered by
       Pedersen \cite{Jesper} and Guo and Shepp \cite{GuoShepp} among others,
       with the aim of computing rational values of perpetual American lookback (Russian) options.
       More recently, Guo and Zervos \cite{GuoZer} derived solutions for discounted optimal
       stopping problems related to the pricing of perpetual American options with certain payoff functions
       depending on the running values of both the initial diffusion process and its associated maximum.
       Glover, Hulley, and Peskir \cite{GHP} provided solutions to optimal stopping problems for integrals
       of functions depending on the running values of both the initial diffusion process and its associated minimum.
              The main feature of the resulting optimal stopping problems
             is that the normal-reflection condition holds for the value function at the diagonal of the state space
       of the two-dimensional continuous Markov process having the initial process and the running extremum
       as its components, which implies the characterisation of the optimal boundaries as the extremal solutions
       of one-dimensional first-order nonlinear ordinary differential equations.

       Asmussen, Avram, and Pistorius \cite{AAP} considered perpetual American options with payoffs depending
       on the running maximum of some L\'evy processes with two-sided jumps having phase-type distributions
       in both directions. Avram, Kyprianou, and Pistorius \cite{AKP} studied exit problems for spectrally
       negative L\'evy processes and applied the results to solving optimal stopping problems for payoff
       functions depending on the running values of the initial processes or their associated maxima.
       Optimal stopping games with payoff functions of such type were considered by Baurdoux and Kyprianou
       \cite{BK3} and Baurdoux, Kyprianou, and Pardo \cite{BKP} within the framework of models based
       on spectrally negative L\'evy processes.
       Other complicated optimal stopping problems for the running maxima were considered in \cite{Gap}
       for a jump-diffusion model with compound Poisson processes with exponentially distributed jumps
       and by Ott \cite{Ott} and Kyprianou and Ott \cite{KyprOtt} (see also Ott \cite{Ott_thesis})
       for a model based on spectrally negative L\'evy processes.
       More recently, Peskir \cite{Pe5a}-\cite{Pe5b} studied optimal stopping problems for three-dimensional Markov
       processes having the initial diffusion process as well as its maximum and minimum as the state space components.
       It was shown that the optimal boundary surfaces depending on the maximum and minimum of the initial process provide
       the maximal and minimal solutions of the associated systems of first-order non-linear partial differential equations.

In this paper, we obtain closed-form solutions to the problems of rational valuation of the perpetual American options on the maximum of the market depth with fixed and floating strikes in an extension of the Black-Merton-Scholes model with path-dependent coefficients.
Such options represent protections for the holders of particularly risky assets, the prices of which can fall deeply, after achieving their historic maxima.
These contracts should be exercised when the maximum drawdown of the underlying asset price rises above the difference of the running maximum and either a certain fixed value or the current value of a certain number of assets.
The closed-form expressions for the rational prices of perpetual American standard put and
call options on the underlying assets in this model were recently computed in \cite{GR3}.
The maximum drawdown process represents the maximum of the difference between the running values of the underlying asset price and its maximum and can therefore be interpreted as the maximum of the market depth. We assume that the price dynamics of the underlying asset are described by a geometric diffusion-type process $X$ with local drift and diffusion coefficients, which essentially depend on the running values of the maximum process $S$ and the maximum drawdown process $Y$.
Such dependence of the dividend and volatility rates on the past dynamics of the asset in the financial market is often used in practice, although it has not been well captured by local or stochastic dividend and volatility models studied in the literature.

It is shown that the optimal exercise times for these options are the first times at which the process $X$ hits certain boundaries depending on the running values of $S$ and $Y$. We derive closed-form expressions for the value functions as solutions of the equivalent free-boundary problems and apply the maximality principle from \cite{Pmax} to describe the optimal boundary surfaces as the {\it minimal} solutions of first-order nonlinear ordinary differential equations.
The starting conditions for these surfaces at the edges of the three-dimensional state space of $(X, S, Y)$ are specified from the solutions of the corresponding optimal stopping problems in the model with the coefficients of the process $X$ depending only on the running maximum process $S$. We also present an explicit solution of the ordinary differential equation corresponding to the case of options with floating strikes in the latter particular model.

The Laplace transforms of the drawdown process and other related characteristics associated with certain classes of the initial processes such as some diffusion models and spectrally positive and negative L\'evy processes were studied by Pospisil, Vecer, and Hadjiliadis \cite{PVH} and by Mijatovi\'c and Pistorius \cite{MP}, respectively.
Diffusion-type processes with given joint laws for the terminal level and supremum at an independent exponential time were
constructed in Forde \cite{Forde}, by allowing the diffusion coefficient to depend on the running values of the initial process
and its running minimum. Other important characteristics for such diffusion-type processes were recently derived by Forde, Pogudin, and Zhang \cite{FPZ}.

      The paper is organized as follows.
       In Section 2, we formulate the associated optimal stopping problems for a necessarily three-dimensional continuous Markov process, which has the underlying asset price and the running values of its maximum and maximum drawdown as the state space components.
       The resulting optimal stopping problems are reduced to their equivalent free-boundary problems for the value functions which satisfy the smooth-fit conditions at the stopping boundaries and the normal-reflection conditions at the edges of the state
       space of the three-dimensional process.
       In Section 3, we obtain closed-form solutions of the associated free-boundary problems and derive first-order nonlinear
       ordinary differential equations for the candidate stopping boundaries. We specify the starting conditions for the
       latter and provide a recursive algorithm to determine the value functions and the optimal boundaries along with their
       intersection lines with the edges of the three-dimensional state space.
       In Section 4, by applying the change-of-variable formula with local time on surfaces from Peskir \cite{Pe1a}, we verify that the resulting solutions of the free-boundary problems provide the expressions for the value functions and the optimal stopping boundaries for the underlying asset price process in the initial problems. Applying an extension of the maximality principle
       from \cite{Pmax} to the three-dimensional optimal stopping problems, we show that the optimal stopping boundaries provide the minimal solutions of the associated ordinary differential equations (see also \cite{Pe5b} for another three-dimensional case).
       The main results of the paper are stated in Theorem \ref{thmSY}.

     \section{\dot Formulation of the problem}

     In this section, we introduce the setting and notation of the three-dimensional optimal stopping problems
     which are related to the pricing of certain perpetual American options and formulate the equivalent free-boundary problems.

     \vspace{6pt}

     {\bf 2.1. Formulation of the problem.}
     For a precise formulation of the problem,
     let us consider a probability space $(\Omega, {\cal F}, P)$
     with a standard Brownian motion $B=(B_t)_{t \ge 0}$ and its natural filtration $({\cal F}_t)_{t \ge 0}$.
     Assume that there exists a continuous process $X=(X_t)_{t \ge 0}$ solving the stochastic differential equation
       \begin{equation}
       \label{dX4}
       dX_t = (\rho-\delta(S_t, Y_t)) \, X_t \, dt + \sigma(S_t, Y_t) \, X_t \, dB_t
       \quad (X_0=x)
       \end{equation}
     where $r > 0$ is a given constant, $\delta(s, y), \sigma(s, y) > 0$
    are some continuously differentiable and bounded functions on $[0, \infty]^2$,
    and $x > 0$ is fixed. Here, the associated with $X$
     {\it running maximum} process $S=(S_t)_{t \ge 0}$ and the corresponding
     {\it running maximum drawdown} process $Y=(Y_t)_{t \ge 0}$ are defined by
        \begin{equation}
        \label{S4}
        S_t = s \vee \max_{0 \le u \le t} X_u \quad \text{and} \quad Y_t = y \vee \max_{0 \le u \le t} (S_u - X_u)
        \end{equation}
for arbitrary $0 < s - y \le x \le s$. Observe that, since the functions $\delta(s, y)$ and
$\sigma(s, y)$ are assumed to be bounded on $[0, \infty]^2$, it follows from the result of
\cite[Chapter~IV, Theorem~4.8]{LS} that there exists a (pathwise) unique solution of the
stochastic differential equation in (\ref{dX4}) which admits the representation
       \begin{equation}
       \label{X4}
       X_t = x \, \exp \bigg( \int_0^t \Big( \rho - \delta(S_u, Y_u) -
       \frac{\sigma^2(S_u, Y_u)}{2} \Big) \, du + \int_0^t \sigma(S_u, Y_u) \, dB_u \bigg)
       \end{equation}
       for all $t \ge 0$.
     We further assume that the resulting continuous diffusion-type process $X$ describes the price of a risky asset on a financial market, where $r$ is the riskless interest rate, $\delta(s, y)$ is the dividend rate paid to the asset holders,
     and $\sigma(s, y)$ is the volatility rate.

     The main purpose of the present paper is
     to derive a closed-form solution to the optimal stopping
     problem for the continuous time-homogeneous (strong) Markov
     process $(X, S, Y)=(X_t, S_t, Y_t)_{t \ge 0}$ given by
        \begin{equation}
        \label{V5b}
        V_*(x, s, y) = \sup_{\tau} E_{x, s, y}
        \big[ e^{- \rho \tau} \, G(X_{\tau}, S_{\tau}, Y_{\tau}) \big]
        \end{equation}
       for any $(x, s, y) \in E^3$, where the supremum is taken over all stopping times $\tau$
       with respect to the natural filtration of $X$, and the payoff function is either
       $G(x, s, y) = (K - s + y)^+$ or $G(x, s, y) = (K x - s + y)^+$.
              Here $E_{x, s, y}$ denotes the expectation
       with respect to the (unique) martingale measure (see, e.g. \cite[Chapter~VII, Section~3g]{FM}),
       under the assumption that the (three-dimensional) process $(X, S, Y)$ defined in
       (\ref{dX4})-(\ref{S4}) starts at $(x, s, y) \in E^3$,
       and $E^3 = \{ (x, s, y) \in \RR^3 \, | \, 0 < s-y \leq x \leq s \}$
       provides the state space for the process $(X, S, Y)$.
              The value of (\ref{V5b}) is then actually a {\it rational} (or {\it no-arbitrage})
       price of a perpetual American option on the maximum of the market depth with payoff
       function either $G(x, s, y) = (K - s + y)^+$ or $G(x, s, y) = (K x - s + y)^+$,
       which corresponds to either the case of fixed strike $K > 0$ or floating strike $K x > 0$, respectively.
       The values of perpetual American standard options were computed in \cite{GR3}
       in the same diffusion-type model.
       The case of perpetual American lookback options with fixed and floating strikes with
       payoff functions $G(x, s, y) = (s - K)^+$ and $G(x, s, y) = (s - K x)^+$
       in the diffusion model of (\ref{dX4})
       with constant coefficients
       $\delta(s, y) = \delta$ and $\sigma(s, y) = \sigma$
       was considered in \cite{Jesper} and \cite{GuoShepp},
       and more complicated $\pi$-options were studied in \cite{GuoZer}.

              \vspace{6pt}

{\bf 2.2. The structure of the optimal stopping times.}
It follows from the general theory of optimal stopping problems
for Markov processes (see, e.g. \cite[Chapter~I, Section~2.2]{PSbook}) that the
optimal stopping time in the problem of (\ref{V5b}) is given by
\begin{equation}
\label{tau0} \tau_* = \inf \{t \ge 0 \, | \, V_*(X_t, S_t, Y_t) = G(X_t, S_t, Y_t) \}
\end{equation}
for the payoff function being either $G(x,s,y) = (K - s + y)^+$ or
$G(x,s,y) = (Kx - s +y)^+$.
The following assertion specifies the structure of the optimal stopping time
$\tau_*$ in (\ref{tau0}) in the both cases of payoff functions.

\begin{lemma} \label{lem}
Suppose that $\delta(s,y), \sigma(s,y) > 0$ are continuously differentiable
bounded functions on $[0,\infty]^2$, and $r>0$ in (\ref{dX4}).
Then, in the optimal stopping problem of (\ref{V5b}) with the payoff function being
either $G(x,s,y) = (K-s+y)^+$ or $G(x,s,y) = (Kx-s+y)^+$, the optimal stopping time
from (\ref{tau0}) has the structure
\begin{equation}
\label{tau*}
\tau_* = \inf \{t \ge 0 \, | \, X_t \geq b_*(S_t, Y_t) \}
\end{equation}
for some function $b_*(s,y)$ to be determined, such that:

(i) for $G(x,s,y) = (K-s+y)^+$, we have
\begin{equation}
\label{b*fix}
b_*(s,y) > s - y
\; \text{ such that }
\; s-y < K
\; , \quad \text{and} \quad
b_*(s,y) > s \geq K + y
\end{equation}
for all $0<y<s$;

(ii) for $G(x,s,y) = (Kx-s+y)^+$, we have
\begin{equation}
\label{b*float}
b_*(s, y) \geq {\underline b}(s, y) \vee (s-y) \vee \frac{s-y}{K}
\; \text{ with } \;
{\underline b}(s, y) = \frac{\rho (s-y)}{\delta(s, y) K}
\end{equation}
for all $0<y<s$.
\end{lemma}

{\bf Proof.} {\bf (i) The case of fixed strike.}
In the case of $G(x, s, y) = (K - s + y)^+$,
it is obvious from the structure of the payoff 
that it is never optimal to stop when $S_t - Y_t \ge K$, for any $t \ge 0$.
In other words,
the set
\begin{equation}
\label{C41}
C' = \{ (x, s, y) \in E^3 \, | \, 0 < K \le s-y \leq x \leq s \} \\
\end{equation}
belongs to the continuation region
\begin{equation}
\label{C40}
C_* = \{ (x, s, y) \in E^3 \, | \, V_*(x, s, y) > (K - s + y)^+ \}.
\end{equation}
It is seen from the solution below that $V_*(x, s, y)$ is continuous,
so that $C_*$ is open.
Then, applying the change-of-variable formula
from \cite{Pe1a} to the function $e^{- r t} (K - s + y)^+$, we get
\begin{align}
\label{GIto4}
&e^{- \rho t} \, (K - S_t + Y_t)^+ = (K - s + y)^+
+ \int_0^t e^{- \rho u} \, I(K > S_u - Y_u) \, dY_u \\
\notag
& - \int_0^t e^{- \rho u} \, I(K > S_u - Y_u) \, dS_u
- \int_0^t e^{- \rho u} \, \rho (K - S_u + Y_u) \, I(K > S_u - Y_u) \, du
\end{align}
where $I(\cdot)$ denotes the indicator function.
It thus follows from the expression in (\ref{GIto4}) that
\begin{align}
\label{Vsigma4}
&E_{x, s, y} \big[ e^{- \rho \tau} \, (K - S_{\tau} + Y_{\tau})^+ \big]
= (K - s + y)^+
+ E_{x, s, y} \bigg[ \int_0^{\tau} e^{- \rho u} \, I(K > S_u - Y_u) \, dY_u \\
\notag
& - \int_0^{\tau} e^{- \rho u} \, I(K > S_u - Y_u) \, dS_u
- \int_0^{\tau} e^{- \rho u} \, \rho (K - S_u + Y_u) \, I(K > S_u - Y_u) \, du  \bigg]
\end{align}
holds for any stopping time $\tau$ and all $0 < s-y \le x \le s$.
Observe that the process $S$ can increase only at the plane
$d_1 = \{(x, s, y) \in \RR^3 \, | \, 0 < x = s \}$,
while the process $Y$ can increase only at the plane
$d_2 = \{(x, s, y) \in \RR^3 \, | \, 0 < x = s - y \}$.
This fact yields through
(\ref{Vsigma4}) that it is never optimal to stop when
$X_t = S_t - Y_t < K$ for $t \ge 0$, so that
the plane $\{ (x, s, y) \in \RR^3 \, | \, 0 < x = s - y < K \}$
belongs to the continuation region in (\ref{C40}).

Since the process $(X, S, Y)$ stays at the same level under the second and
third coordinates, as long as it fluctuates between the planes $d_1$ and $d_2$,
it is clear that we should not let the process $X$ grow up too
much, since it might hit the plane
$d_1$, that happens when $X_t = S_t$ for $t \ge 0$,
and thereby increase its running maximum $S$ and thus decrease the payoff $(K - s + y)^+$.
Moreover, it is not optimal to let the process $X$ run too much away from the plane
$d_2$, since the cost of waiting until it comes
back to the plane $d_2$, that happens when $X_t = S_t - Y_t$ for $t \ge 0$,
and thereby increase its running maximum drawdown $Y$ and thus
the payoff $(K - s + y)^+$ could be too large due to the negative last term in (\ref{Vsigma4})
and the discounting factor in (\ref{V5b}).
It follows from the structure of the value function in (\ref{V5b}) with the processes of
(\ref{dX4})-(\ref{S4}) and the convex payoff function $G(x, s, y)$ that $V_*(x, s, y)$ is convex in $x$
on the interval $(s-y, s)$. Then, standard geometric arguments imply that
there exists a function $b_*(s, y)$ such that the continuation region in (\ref{C40})
consists of (\ref{C41}) and the set
\begin{equation}
\label{C42}
C'' = \{(x, s, y) \in E^3 \, | \, s - y \leq x < b_*(s, y) \;\; \text{and} \;\; s - y < K \}
\end{equation}
while the corresponding stopping region is the closure of the set
\begin{equation}
\label{D4}
D_* = \{(x, s, y) \in E^3 \, | \, b_*(s, y) < x \leq s \;\; \text{and} \;\; s - y < K \}
\end{equation}
with $b_*(s,y) > s - y$ for all $0 < y < s$ such that $0 < s - y < K$.

\vspace{3pt}


{\bf (ii) The case of floating strike.}
In the case of $G(x, s, y) = (K x - s + y)^+$,
it is obvious from the structure of the payoff
that it is never optimal to stop when $K X_t \le S_t - Y_t$,
for any $t \ge 0$. In other words, this fact shows that
the set
\begin{equation}
\label{C51}
C' = \{ (x, s, y) \in E^3 \, | \, 0 < K x \le s - y \}
\end{equation}
which exists for all $0 < K \leq 1$, belongs to the continuation region
\begin{equation}
\label{C50}
C_* = \{ (x, s, y) \in E^3 \, | \, V_*(x, s, y) > (K x - s + y)^+ \}.
\end{equation}
It is seen from the solution below that $V_*(x, s, y)$ is continuous,
so that $C_*$ is open.
Then, applying the change-of-variable formula from \cite{Pe1a}
to the function $e^{- r t} (K x - s + y)^+$, we get
\begin{align}
\label{GIto5}
&e^{- \rho t} \, (K X_t - S_t + Y_t)^+ = (K x - s + y)^+ + N_t \\
\notag
&+ \int_0^t e^{- \rho u} \, H(X_u, S_u, Y_u) \, I(K X_u > S_u - Y_u) \, du
- \int_0^t e^{- \rho u} \, I(K X_u > S_u - Y_u) \, dS_u  \\
\notag
&+ \int_0^t e^{- \rho u} \, I(K X_u > S_u - Y_u) \, dY_u +
\frac{1}{2} \int_0^t e^{- \rho u} \, K \, I(K X_u = S_u - Y_u) \, d \ell^K_u(X)
\notag
\end{align}
where we set $H(x, s, y) = \rho (s-y) - K \delta(s, y) x$, and
the process $\ell^K(X) = (\ell^K_t(X))_{t \ge 0}$ is the local time of
$X$ at the plane $\{ (x, s, y) \in \RR^3 \, | \, 0 < K x = s - y \}$ given by
\begin{equation}
\label{locc5a}
\ell^K_t(X) = \lim_{\eps \downarrow 0} \frac{1}{2 \eps} \int_0^t
I(-\eps < K X_u - S_u + Y_u < \eps) \, \sigma^2(S_u, Y_u) \, X^2_u \, du
\end{equation}
as a limit in probability.
Here, the process $N = (N_t)_{t \ge 0}$ defined by
\begin{equation}
\label{Mlocc5a}
N_t = \int_0^t e^{- r u} \, I(K X_u > S_u - Y_u) \,
\sigma(S_u, Y_u) \, X_u \, d{B}_u
\end{equation}
is a continuous square integrable martingale under $P_{x, s, y}$.
Hence, applying Doob's optional sampling theorem
(see, e.g. \cite[Chapter~I, Theorem~3.22]{KS}) and using the expression in (\ref{GIto5}),
we get that
\begin{align}
\label{Vsigma5}
&E_{x, s, y} \big[ e^{- \rho \tau} \, (K X_{\tau} - S_{\tau} +
Y_{\tau})^+ \big] = (K x - s + y)^+ \\
\notag
&+ E_{x, s, y} \bigg[ \int_0^{\tau} e^{- \rho u} \, H(X_u, S_u, Y_u) \,
I(K X_u > S_u - Y_u) \, du - \int_0^{\tau} e^{- \rho u} \, I(K X_u > S_u - Y_u) \, dS_u  \\
\notag
&\phantom{+ E_{x, s, y} \bigg[ \;\:}+
\int_0^{\tau} e^{- \rho u} \, I(K X_u > S_u - Y_u) \, dY_u + \frac{1}{2}
\int_0^{\tau} e^{- \rho u} \, K \,
I(K X_u = S_u - Y_u) \, d \ell^K_u(X) \bigg]
\end{align}
holds for any stopping time $\tau$ and all $0 < s-y \le x \le s$.
It is seen from (\ref{Vsigma5}) that it is never optimal to stop when
$K X_t > S_t - Y_t$ and either $H(X_t, S_t, Y_t) > 0$ or $X_t = S_t - Y_t$
holds for $t \ge 0$, where the latter condition is true since the process $Y$ can only
increase at the plane $d_2$.
In other words,
the set
\begin{equation}
\label{C52}
C'' = \{ (x, s, y) \in E^3 \, | \, (s-y) / K < x < {\underline b}(s,y)
\;\; \text{or} \;\; x = s - y > (s-y) / K \}
\end{equation}
with ${\underline b}(s,y) = {\rho (s-y)}/({\delta(s, y) K})$
for $0 < y < s$, belongs to the continuation region in (\ref{C50}).
Note that, the set in (\ref{C52}) exists only if
$K>1$ holds or if
$0 < K \leq 1$ and $\delta(s,y) < r$ holds.

Let us now fix some $(x, s, y) \in C_*$ from the continuation region in
(\ref{C50}) and let $\tau_* = \tau_*(x, s, y)$ denote the optimal stopping
time in the problem of (\ref{V5b}) with $G(x, s, y) = (K x - s + y)^+$.
Then, by means of the results of general optimal stopping theory for Markov
processes (see, e.g. \cite[Chapter~I, Section~2.2]{PSbook}),
we conclude from the expression in (\ref{Vsigma5}) that
\begin{align}
\label{Vsigma5c}
&V_*(x, s, y) - (K x - s + y)^+ \\
\notag
&= E_{x, s, y} \bigg[ \int_0^{\tau_*} e^{- \rho u} \, H(X_u, S_u, Y_u) \,
I(K X_u > S_u - Y_u) \, du - \int_0^{\tau_*} e^{- \rho u} \, I(K X_u > S_u - Y_u) \, dS_u \\
\notag
&\phantom{= E_{x, s, y} \bigg[ \;\:}+ \int_0^{\tau_*} e^{- \rho u} \, I(K X_u > S_u - Y_u)
\, dY_u + \frac{1}{2} \int_0^{\tau_*} e^{- \rho u} \, K \, I(K X_u = S_u - Y_u) \, d \ell^K_u(X) \bigg] > 0
\end{align}
holds. Hence, taking any $x'$ such that ${\underline b}(s, y) \vee (s - y)
< x' < x$ and using the explicit expression for the process $X$ through its
starting point in (\ref{X4}) and also the structure of the maximum drawdown
process $Y$ in (\ref{S4}), we obtain from (\ref{Vsigma5}) that the inequalities
\begin{align}
\label{Vsigma5cx'}
&V_*(x', s, y) - (K x' - s + y)^+ \\
\notag
&\geq E_{x', s, y} \bigg[ \int_0^{\tau_*} e^{- \rho u} \, H(X_u, S_u, Y_u) \,
I(K X_u > S_u - Y_u) \, du - \int_0^{\tau_*} \, e^{- \rho u} \, I(K X_u > S_u - Y_u) \,
dS_u \\
\notag
&\phantom{= E_{x', s, y} \bigg[ \;\:}+ \int_0^{\tau_*} e^{- \rho u} \, I(K X_u > S_u - Y_u)
\, dY_u + \frac{1}{2} \int_0^{\tau_*} e^{- \rho u} \, K \, I(K X_u > S_u - Y_u) \, d \ell^K_u(X) \bigg] \\
\notag
&\ge E_{x, s, y} \bigg[ \int_0^{\tau_*} e^{- \rho u} \, H(X_u, S_u, Y_u) \,
I(K X_u > S_u - Y_u) \, du - \int_0^{\tau_*} e^{- \rho u} \, I(K X_u > S_u - Y_u) \,
dS_u \\
\notag
&\phantom{= E_{x, s, y} \bigg[ \;\:}+ \int_0^{\tau_*} e^{- \rho u} \, I(K X_u > S_u - Y_u)
\, dY_u + \frac{1}{2} \int_0^{\tau_*} e^{- \rho u} \, K \, I(K X_u = S_u - Y_u) \, d \ell^K_u(X) \bigg]
\end{align}
are satisfied.
Thus,
by virtue of the inequality in (\ref{Vsigma5c}),
we see that $(x', s, y) \in C_*$.
Taking into account the convexity of the function $G(x, s, y)$ and thus of $V_*(x, s, y)$
in $x$ on the interval $(s-y, s)$, we may therefore conclude from the standard geometric
arguments that
there exists a function $b_*(s, y)$ such that
the continuation region in (\ref{C50}), which
consists of the regions in (\ref{C51}) and (\ref{C52}), takes the form
\begin{equation}
\label{C5}
C_* = \{(x, s, y) \in E^3 \, | \, s-y \le x < b_*(s,y) \}
\end{equation}
while the corresponding stopping region is the closure of the set
\begin{equation}
\label{D5}
D_* = \{(x, s, y) \in E^3 \, | \, b_*(s,y) < x \leq s \}
\end{equation}
with
$b_*(s, y) \geq {\underline b}(s, y) \vee (s-y) \vee ((s-y)/K)$
and
${\underline b}(s, y) = {\rho (s-y)}/({\delta(s, y) K})$ for $0 < y < s$. $\square$

        \vspace{6pt}

{\bf 2.3. The free-boundary problem.}
By means of standard arguments based on the application of
It\^o's formula, it is shown that the infinitesimal operator $\LL$ of the
process $(X, S, Y)$ acts on a function $F(x, s, y)$ from the class
$C^{2,1,1}$ on the interior of $E^3$ according to the rule
     \begin{equation}
\label{LF}
(\LL F)(x, s, y) = (\rho - \delta(s, y)) \, x \, \partial_x F(x, s, y)
+ \frac{\sigma^2(s, y)}{2} \, x^2 \, \partial^2_{xx} F(x, s, y)
\end{equation}
for all $0 < s-y < x < s$.
      In order to find analytic expressions for the unknown
      value function $V_*(x, s, y)$ from (\ref{V5b}) and the unknown
      boundary $b_*(s, y)$ from (\ref{tau*}),
      let us build on the results of general theory of optimal stopping
      problems for Markov processes (see, e.g. 
      \cite[Chapter~IV, Section~8]{PSbook}). We can reduce the optimal
      stopping problem of (\ref{V5b}) to the equivalent free-boundary problem
      for $V_*(x, s, y)$ and $b_*(s, y)$ given by
       \begin{align}
       \label{LV}
       &({\LL} V)(x, s, y) = \rho V(x, s, y)
       \quad \text{for} \quad (x, s, y) \in C \\
       \label{CF}
       &V(x, s, y) \big|_{x = b(s, y)-} = G(b(s, y), s, y) \\
       \label{VD}
       &V(x, s, y) = G(x, s, y) \quad \text{for} \quad (x, s, y) \in D \\
       \label{VC}
       &V(x, s, y) > G(x, s, y) \quad \text{for} \quad (x, s, y) \in C \\
       \label{LVD}
       &({\LL} V)(x, s, y) < \rho V(x, s, y) \quad \text{for} \quad (x, s, y) \in D
      \intertext{where $C$ is defined as $C' \cup C''$ in (\ref{C41}) and (\ref{C42}) or $C_*$ in (\ref{C5}),
and $D$ is defined as $D_*$ in (\ref{D4}) or (\ref{D5}),
for the payoff function $G(x,s,y) = (K-s+y)^+$ or $G(x,s,y) = (Kx-s+y)^+$, respectively,
      with $b(s, y)$ instead of $b_*(s, y)$.
      The instantaneous-stopping condition in (\ref{CF}) is satisfied, when $s - y \le b(s, y) \le s$
      holds, for each $0 < y < s$. Observe that the superharmonic characterisation of the value function
      (see \cite{Dyn} 
      and \cite[Chapter~IV, Section~9]{PSbook}) implies
      that $V_*(x, s, y)$ is the smallest function satisfying
      (\ref{LV})-(\ref{VC}), with the boundary $b_*(s, y)$.
      Moreover, we further assume that the normal-reflection and the smooth-fit conditions}
       \label{NRY}
       &\partial_y V(x, s, y) \big|_{x = (s-y)+} = 0 \quad \text{and} \quad
       \partial_x V(x, s, y) \big|_{x = b(s, y)-} = \partial_x G(b(s, y), s, y)
       \intertext{are satisfied, when $s - y < b(s, y) < s$ holds, for each
$0 < y < s$.
       Otherwise, we assume that the normal-reflection conditions}
\label{NRboth}
&\partial_y V(x, s, y) \big|_{x = (s-y)+} = 0 \quad \text{and} \quad
\partial_s V(x, s, y) \big|_{x = s-} = 0
\end{align}
are satisfied, when $b(s, y) > s$ holds, for each $0 < y < s$.

Observe that when the inequalities
$S_t - Y_t < X_t < b(S_t, Y_t) < S_t$ are satisfied for $t \ge 0$,
the process $X$ can increase towards the boundary $b(S_t, Y_t)$ in a
continuous way, so that we can assume that the smooth-fit condition of
(\ref{NRY}) holds for the candidate value function $V(x, s, y)$ at the
boundary $b(s, y)$.
Such assumptions are naturally applied to determine the solutions of the
free-boundary problems, which provide the solutions of associated optimal
stopping problems
(see \cite[Chapter~IV, Section~9]{PSbook} for an explanation and proofs).
On the other hand, when either the inequalities
$S_t - Y_t < X_t < b(S_t, Y_t) < S_t$
or $S_t - Y_t < X_t < S_t < b(S_t, Y_t)$ are satisfied for $t \ge 0$,
the process $X$ can increase or decrease towards the planes $d_1$ or $d_2$,
respectively, in a continuous way. In this case, it follows
from the property of the infinitesimal operator
of the process $(X, S, Y)$ that
the normal-reflection conditions of (\ref{NRY}) and (\ref{NRboth}) hold for the candidate
value function $V(x, s, y)$ at the planes $d_1$ and $d_2$.
These conditions are used to derive first-order ordinary differential equations for the
candidate boundaries of the corresponding optimal stopping problems
(see \cite{DSS}, \cite{GP1}-\cite{GP2},
\cite{SS1}-\cite{SS2}, \cite{Jesper}, \cite{GuoShepp}, \cite{GuoZer}, and \cite{GHP} among others,
and \cite{Pmax} and \cite[Chapter~IV, Section~13]{PSbook} for an explanation and further references).
We follow the classical approach and apply the smooth-fit and normal-reflection conditions
from (\ref{NRY}) and (\ref{NRboth}) to find closed-form expressions for the candidate
value functions as well as ordinary differential equations for the boundaries, and then
verify in Theorem \ref{thmSY} below that
the obtained solutions to the free-boundary problem provide the value
functions and the optimal stopping boundaries in the original problems.


     \vspace{6pt}

{\bf 2.4. Some remarks.}
Let us finally note some facts about the value function $V_*(x, s, y)$ in (\ref{V5b}),
which will then be used in order to specify the asymptotic behavior of the
boundary $b_*(s, y)$ from (\ref{tau*}).

\vspace{3pt}

{\bf (i) The case of fixed strike.} In the case of $G(x, s, y) = (K - s + y)^+$,
we observe from (\ref{V5b}) that the inequalities
\begin{equation}
\label{WVW3}
0 \le (K - s + y)^+
\le \sup_{\tau} E_{x, s, y} \big[ e^{- \rho \tau} \, (K - S_{\tau} + Y_{\tau})^+ \big]
\le K
\end{equation}
hold for all $(x, s, y) \in E^3$.
Thus, setting $x = s - y$ into (\ref{WVW3}) and letting $y$ increase to $s$, we get that the property
\begin{equation}
\label{WVW5}
\liminf_{y \uparrow s} V_*(s-y, s, y) = \limsup_{y \uparrow s} V_*(s-y, s, y) = K
\end{equation}
holds for all $s > 0$.


\vspace{3pt}

{\bf (ii) The case of floating strike.}
In the case of $G(x, s, y) = (K x - s + y)^+$, we observe
from (\ref{V5b}) that the inequalities
\begin{align}
\label{WVW6}
0 \le (K \, x - s + y)^{+}
\le \sup_{\tau} E_{x, s, y} \big[ e^{- \rho \tau} \, (K \, X_{\tau} - S_{\tau} +
Y_{\tau})^+ \big]
\le K \, \sup_{\tau} E_{x, s, y} \big[ e^{- \rho \tau} \, X_{\tau} \big]
\end{align}
imply that the expressions
\begin{equation}
\label{WVW7}
0 \le (K \, x - s + y)^{+} \le V_*(x, s, y) \le K \, x
\end{equation}
hold for all $(x, s, y) \in E^3$, where
the third inequality in (\ref{WVW7}) follows from the optimal immediate stopping
of the problem in the right-hand side of (\ref{WVW6}).
Thus, setting $x = s - y$ in (\ref{WVW7}) and letting $y$ increase to $s$, we obtain that the property
\begin{equation}
\label{WVW9}
0 \le (K - 1)^{+} \le \lim_{y \uparrow s} \frac{V_*(s-y, s, y)}{s-y} \le K
\end{equation}
holds for all $s > 0$.

       \section{\dot Solution of the free-boundary problem}

In this section, we obtain closed-form expressions for the value functions
$V_*(x, s, y)$ in (\ref{V5b}) associated with the options on the maximum of the
market depth with fixed and floating strikes, and derive first-order nonlinear
ordinary differential equations for the optimal exercise boundaries $b_*(s, y)$ from
(\ref{tau*}), as solutions to the free-boundary problem in (\ref{LV})-(\ref{NRboth}).
The analysis performed in this section also provides a recursive algorithm to determine
the candidate value functions and optimal stopping boundaries as well as their
intersection lines with the edges of the three-dimensional state space.

\vspace{6pt}

{\bf 3.1. The general solution of the ordinary differential equation.}
We first observe that the general solution of the equation in (\ref{LV})
has the form
\begin{equation}
\label{V0}
V(x, s, y)
= C_1(s, y) \, x^{\gamma_1(s, y)} + C_2(s, y) \, x^{\gamma_2(s, y)}
\end{equation}
where $C_i(s, y)$, $i = 1, 2$, are some arbitrary
continuously differentiable functions and $\gamma_i(s, y)$,
$i = 1, 2$, are given by
\begin{equation}
\label{gamma12}
\gamma_{i}(s, y) = \frac{1}{2} - \frac{\rho - \delta(s, y)}{\sigma^2(s, y)}
- (-1)^i \sqrt{ \bigg( \frac{1}{2} - \frac{\rho - \delta(s, y)}
{\sigma^2(s, y)} \bigg)^2 + \frac{2\rho}{\sigma^2(s, y)}}
\end{equation}
so that $\gamma_2(s, y) < 0 < 1 < \gamma_1(s, y)$ holds for all $0 < y < s$.
Hence, applying the instantaneous-stopping condition from (\ref{CF})
to the function in (\ref{V0}), we get that the equality
\begin{equation}
\label{B31b}
C_1(s, y) \, b^{\gamma_1(s, y)}(s, y)
+ C_2(s, y) \, b^{\gamma_2(s, y)}(s, y) = G(b(s, y), s, y)
\end{equation}
is satisfied, when $s - y \le b(s, y) \le s$ holds, for each $0 < y < s$. Moreover, using
the smooth-fit condition from the right-hand part of (\ref{NRY}),
we obtain that the equality
\begin{equation}
\label{B31bb}
C_1(s, y) \, \gamma_1(s, y) \, b^{\gamma_1(s, y)}(s, y) + C_2(s, y) \,
\gamma_2(s, y) \, b^{\gamma_2(s, y)}(s, y) = \partial_x G(b(s, y), s, y) \, b(s, y)
\end{equation}
is satisfied, when $s - y < b(s, y) < s$ holds, for each $0 < y < s$.
Finally, applying the normal-reflection conditions
from (\ref{NRboth}) to the function in (\ref{V0}), we obtain that the equalities
\begin{align}
\label{B31c}
&\sum_{i=1}^{2} \Big( \partial_s C_i(s, y) \, s^{\gamma_i(s,y)} + C_i(s,y)
\, \partial_s \gamma_i(s,y) \, s^{\gamma_i(s,y)} \ln s \Big) = 0 \\
\label{B31d}
&\sum_{i=1}^{2} \Big( \partial_y C_i(s, y) \, (s-y)^{\gamma_i(s,y)} +
C_i(s, y) \, \partial_y \gamma_i(s,y) \, (s-y)^{\gamma_i(s,y)} \ln (s-y) \Big) = 0
\end{align}
are satisfied, when $s < b(s, y)$ and $s-y < b(s, y)$ holds, respectively, for each $0 < y < s$.
Here, the partial derivatives $\partial_s \gamma_{i}(s, y)$ and $\partial_y \gamma_{i}(s, y)$
take the form
\begin{align}
\label{gammas}
&\partial_s \gamma_{i}(s, y) = \phi(s, y) - (-1)^i \frac{
\phi(s, y) (\gamma_1(s, y) - \gamma_2(s, y)) \sigma^3(s, y) - 2 \rho \partial_s \sigma(s, y)}
{\sigma^2(s, y) \sqrt{(\gamma_1(s, y) - \gamma_2(s, y))^2 \sigma^2(s, y) + 2 \rho}} \\
\label{gammay}
&\partial_y \gamma_{i}(s, y) = \psi(s, y) - (-1)^i \frac{
\psi(s,y) (\gamma_1(s, y) - \gamma_2(s, y)) \sigma^3(s, y) - 2 \rho \partial_y \sigma(s, y)}
{\sigma^2(s, y) \sqrt{(\gamma_1(s, y) - \gamma_2(s, y))^2 \sigma^2(s, y) + 2 \rho}}
\end{align}
for $i = 1, 2$, and the functions $\phi(s, y)$ and $\psi(s, y)$ are defined by
\begin{align}
\label{phi}
&\phi(s,y) = \frac{\sigma(s, y) \partial_s \delta(s, y) + 2
(\rho - \delta(s, y)) \partial_s \sigma(s, y)}
{\sigma^3(s, y)} \\
\label{psi}
&\psi(s,y) = \frac{\sigma(s, y) \partial_y \delta(s, y) + 2
(\rho - \delta(s, y)) \partial_y \sigma(s, y)}
{\sigma^3(s, y)}
\end{align}
for $0 < y < s$.

\vspace{6pt}


{\bf 3.2. The solution to the problem in the $\delta(s)$ and $\sigma(s)$-setting.}
We begin with the case in which $\delta(s, y) = \delta(s)$ and
$\sigma(s, y) = \sigma(s)$ holds in (\ref{dX4}), and thus, we can define
the functions $\beta_i(s) = \gamma_i(s, y)$, $i = 1, 2$, as in (\ref{gamma12}).
Then, the general solution $V(x, s, y)$ of the equation in (\ref{LV}) has the form of
(\ref{V0}) with $\gamma_i(s, y) = \beta_i(s)$, for $i = 1, 2$.
Recall that the border planes of the state space
$E^3 = \{ (x, s, y) \in \RR^3 \, | \, 0 < s-y \leq x \leq s \}$ are $d_1 = \{(x, s, y)
\in \RR^3 \, | \, 0 < x = s \}$ and $d_2 = \{(x, s, y) \in \RR^3 \, | \, 0 < x =
s-y \}$, as well as that the second and third components of the process $(X, S, Y)$
can increase only at the planes $d_1$ and $d_2$, that is, when $X_t = S_t$ and
$X_t = S_t - Y_t$ for $t \ge 0$, respectively.

\vspace{3pt}

{\bf (i) The case of fixed strike.}
Let us first consider the payoff function $G(x,s,y) = (K-s+y)^+$ in (\ref{V5b}).
In this case, solving the system of equations in (\ref{B31b})-(\ref{B31bb}),
we obtain that the function in (\ref{V0}) admits the representation
\begin{equation}
\label{V34p}
V(x, s, y; b(s, y)) = C_1(s, y; b(s, y)) \, x^{\beta_1(s)} +
C_2(s, y; b(s, y)) \, x^{\beta_2(s)}
\end{equation}
for $0 < s-y \le x < b(s, y) \leq s$, with
\begin{equation}
\label{Ci34p}
C_i(s, y; b(s, y)) = \frac{\beta_{3-i}(s) (K - s + y) }
{(\beta_{3-i}(s) - \beta_i(s)) b(s, y)^{\beta_i(s)}}
\end{equation}
for all $0 < y < s$ and $i = 1, 2$.
Hence, assuming that the boundary function $b(s, y)$ is continuously differentiable,
we apply the condition of (\ref{B31d}) to the functions
$C_i(s, y) = C_i(s, y; b(s, y))$, $i = 1, 2$, in (\ref{Ci34p}) to obtain that the
boundary solves the first-order nonlinear ordinary differential equation
\begin{align}
\label{g'34p}
\partial_y b(s, y) &= \sum_{i=1}^{2} \frac{b(s, y)}{\beta_i(s) (K-s+y)}
\bigg( \frac{((s-y)/b(s, y))^{\beta_i(s)}}{((s-y)/b(s, y))^{\beta_{i}(s)}
- ((s-y)/ b(s, y))^{\beta_{3-i}(s)}} \bigg)
\end{align}
for $0 < y < s$.
Taking into account the condition in (\ref{WVW5}) for the value function in
(\ref{V34p})-(\ref{Ci34p}), we conclude after some straightforward calculations
that $b(s, y) \sim g_*(s) (s-y)$ should hold as $y \uparrow s$, where $g_*(s)$
is the unique solution of the equation
\begin{equation}
\label{gfixp}
\sum_{i = 1}^2 \beta_i(s) \, (g^{- \beta_{3-i}(s)}(s) - 1) = 0
\end{equation}
which is given by $g_*(s) = 1$ for all $s > 0$.
Thus, any candidate solution of the differential equation in (\ref{g'34p}) should satisfy the condition
\begin{equation}
\label{startfixp}
\lim_{y \uparrow s} \frac{b(s, y)}{s-y} = 1
\end{equation}
for all $s > 0$.

For any $s > 0$ fixed, let us now consider a candidate solution $b(s, y)$
of the first-order ordinary differential equation of (\ref{g'34p}),
satisfying the starting condition of (\ref{startfixp}),
given that this solution stays strictly above the plane $d_2$ for all $0 < y < s$
such that $s-y < K$, and strictly above
\begin{picture}(160,108)
\put(20,15){\begin{picture}(120,92)

\put(0,0){\line(1,0){120}} \put(120,0){\line(0,1){92}}
\put(0,0){\line(0,1){92}} \put(0,92){\line(1,0){120}}   

\put(10,10){\vector(1,0){100}}
\put(10,10){\vector(0,1){80}}    

\put(81,9){\line(0,1){72}} 

\put(80,6){$s$}

\put(9,81){\line(1,0){72}} 

\put(10,81){\line(1,-1){71}} 

\put(6,80){${s}$}

\put(60,50){\line(-1,0){19}}
\put(60,49.9){\line(-1,0){19}}
\put(60,50.1){\line(-1,0){19}}

\put(41,50){\line(-1,1){10}}
\put(41,49.9){\line(-1,1){10}}
\put(41,50.1){\line(-1,1){10}}
\put(41,49.8){\line(-1,1){10}}
\put(41,50.2){\line(-1,1){10}}

\put(65,53){\line(-1,0){27}}
\put(65,52.9){\line(-1,0){27}}
\put(65,53.1){\line(-1,0){27}}

\put(50,56.5){\line(-1,0){15.5}}
\put(50,56.4){\line(-1,0){15.5}}
\put(50,56.6){\line(-1,0){15.5}}

\put(81,60){\line(-1,0){50}}
\put(81,59.9){\line(-1,0){50}}
\put(81,60.1){\line(-1,0){50}}

\qbezier[5](81,60)(83,60)(85,60)	
\qbezier[80](85,10)(85,47.5)(85,85)
\qbezier[75](10,85)(47.5,85)(85,85)
\qbezier[95](10,85)(47.5,50)(85,10)

\put(85,9){\line(0,1){2}} \put(84,6){$s'$}
\put(9,85){\line(1,0){2}} \put(6,84){$s'$}

\put(29,20){\vector(1,1){20.5}} \put(19,17){$x = s-y$}

\put(53,50){${\bf .}$}
\put(53,49.9){${\bf .}$}	\put(53.1,49.9){${\bf .}$}	\put(52.9,49.9){${\bf .}$}
\put(53,49.8){${\bf .}$}	\put(53.1,49.8){${\bf .}$} 	\put(52.9,49.8){${\bf .}$}
\put(53,49.7){${\bf .}$}	\put(53.1,49.7){${\bf .}$}	\put(52.9,49.7){${\bf .}$}
\put(53,49.6){${\bf .}$}

\put(51,46){$(x,s,y)$}

\put(0,0){\qbezier(10,81)(19,70)(81,63)} 
\put(0,0){\qbezier(81,63)(92,61)(81,55)}
\put(0,0){\qbezier(81,55)(51,42)(81,30)}
\put(0,0){\qbezier(81,30)(95,24)(87,10)}
\put(0,0.1){\qbezier(10,81)(19,70)(81,63)} 
\put(0,0.1){\qbezier(81,63)(92,61)(81,55)}
\put(0,0.1){\qbezier(81,55)(51,42)(81,30)}
\put(0,0.1){\qbezier(81,30)(95,24)(87,10)}
\put(0.1,0){\qbezier(10,81)(19,70)(81,63)} 
\put(0.1,0){\qbezier(81,63)(92,61)(81,55)}
\put(0.1,0){\qbezier(81,55)(51,42)(81,30)}
\put(0.1,0){\qbezier(81,30)(95,24)(87,10)}

\qbezier[71](10,29.9)(47.5,29.9)(81,29.9)
\qbezier[71](10,63)(45.5,63)(81,63)      
\qbezier[71](10,55)(45.5,55)(81,55)

\put(.5,61.4){${\tilde y}_{1}(s)$}
\put(.5,53.4){${\tilde y}_{2}(s)$}
\put(.5,28.4){${\tilde y}_{3}(s)$}

\put(95,79){\vector(-2,-1){27}} \put(96,78){$b(s,y)$}

\qbezier[20](61,10)(61,20)(61,30)
\put(61,9){\line(0,1){2}} \put(59,5.7){$K$}

\put(6,88){$y$} \put(108,6){$x$}
\end{picture}}
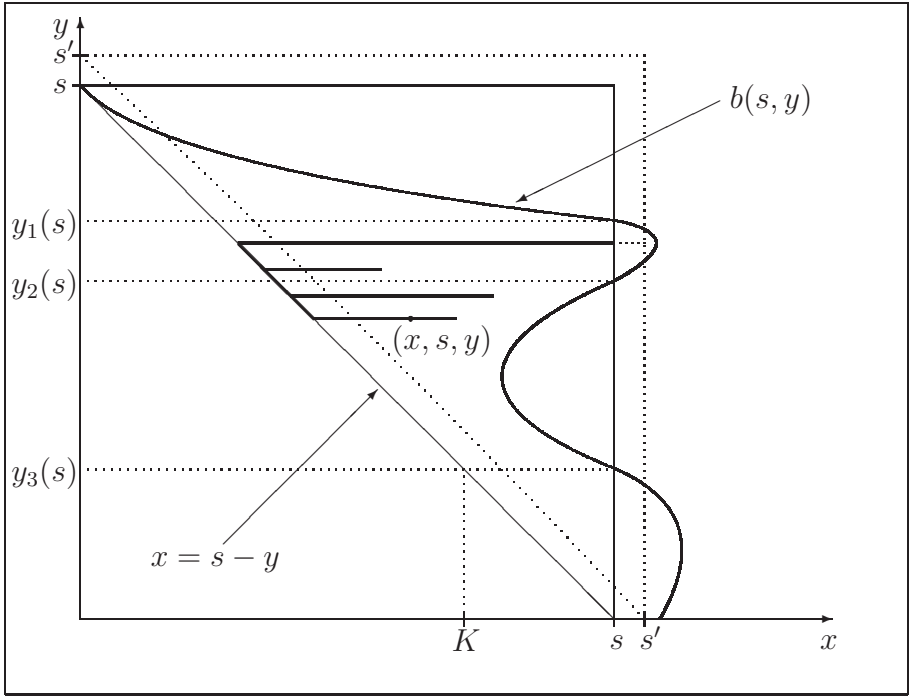
\put(17,10){\small{{\bf Figure 1.} A computer drawing of the
state space of the process $(X,S,Y)$,}} \\
\put(17,5) {\small{for some $s$ fixed, which increases to
$s'$, and the boundary function $b(s,y)$.}}
\end{picture}

\noindent
the plane $d_1$ for all $0 < y < s$ such that $s-y \geq K$.
These assumptions for the boundary function $b(s, y)$ follow from the
structure of the continuation region $C' \cup C''$ in (\ref{C41}) and
(\ref{C42}), which results to the expressions of (\ref{b*fix}) in
Lemma \ref{lem}.
Then, we put ${\tilde y}_0(s) = s$ and define a decreasing sequence
$({\tilde y}_n(s))_{n \in \NN}$ such that the boundary $b(s, y)$ exits
the region $E^3$ from the side of $d_1$ at the points
$(s, s, {\tilde y}_{2l-1}(s))$ and
enters $E^3$ downwards at the points $(s, s, {\tilde y}_{2l}(s))$.
Namely, we define
${\tilde y}_{2l-1}(s) = \sup \{ y < {\tilde y}_{2l-2}(s) \, | \, b(s, y) > s \}$
and ${\tilde y}_{2l}(s) = \sup \{ y < {\tilde y}_{2l-1}(s) \, | \, b(s, y) \le s \}$,
whenever they exist, and put ${\tilde y}_{2l-1}(s) = {\tilde y}_{2l}(s) = 0$,
$l \in \NN$, otherwise.
It follows from (\ref{b*fix}) that
$s-K \le {\tilde y}_{2l-1}(s) < {\tilde y}_{2l-2}(s) \leq s$, for $l = 1,\ldots,{\tilde l}'$,
where ${\tilde l}' = \max \{l \in \NN \,| \, s-{\tilde y}_{2l-1}(s) \le K \}$.
Therefore, the candidate value function admits the expression of (\ref{V34p})-(\ref{Ci34p})
in the regions
\begin{equation}
\label{tilR3_2l-1p}
{\tilde R}_{2l-1} = \{ (x,s,y) \in E^3 \, | \, {\tilde y}_{2l-1}(s) < y \leq {\tilde y}_{2l-2}(s) \}
\end{equation}
for $l = 1, \ldots, {\tilde l}'$ (see Figure 1 above).

On the other hand, the candidate value function
takes the form of (\ref{V0}) with $C_{i}(s, y)$, $i = 1, 2$,
solving the linear system of first-order partial differential equations
in (\ref{B31c})-(\ref{B31d}), in the regions
\begin{equation}
\label{tilR3_2lp}
{\tilde R}_{2l} = \{ (x,s,y) \in E^3 \, | \, {\tilde y}_{2l}(s) < y \leq {\tilde y}_{2l-1}(s) \}
\end{equation}
for $l = 1, \ldots, {\tilde l}'$, which belong to the continuation region $C'\cup C''$
given in (\ref{C41}) and (\ref{C42}).
Note that, the process $(X, S, Y)$ can enter the region ${\tilde R}_{2l}$
in (\ref{tilR3_2lp}) from ${\tilde R}_{2l+1}$ in (\ref{tilR3_2l-1p}), for some
$l = 1,\ldots,{\tilde l}'-1$, only through the point
$(s - {\tilde y}_{2l}(s), s, {\tilde y}_{2l}(s))$ and can exit the region
${\tilde R}_{2l}$ passing to the region ${\tilde R}_{2l-1}$, for some
$l = 1,\ldots,{\tilde l}'$, only through the point $(s - {\tilde y}_{2l-1}(s),
s, {\tilde y}_{2l-1}(s))$, by hitting the plane $d_2$, so that increasing
its third component $Y$.
Thus, the candidate value function should be continuous at the points
$(s - {\tilde y}_{2l}(s), s, {\tilde y}_{2l}(s))$ and $(s - {\tilde y}_{2l-1}(s),
s, {\tilde y}_{2l-1}(s))$, that is expressed by the equalities
\begin{align}
\label{4condcall1}
&C_1(s, {\tilde y}_{2l}(s)+) \, ((s - {\tilde y}_{2l}(s))-)^{\beta_1(s)}
+ C_2(s, {\tilde y}_{2l}(s)+) \, ((s - {\tilde y}_{2l}(s))-)^{\beta_2(s)} \\ \nonumber
&= V(s - {\tilde y}_{2l}(s), s, {\tilde y}_{2l}(s); b(s, {\tilde y}_{2l}(s))) \\
\label{4condcall3}
&C_1(s, {\tilde y}_{2l-1}(s)) \, (s - {\tilde y}_{2l-1}(s))^{\beta_1(s)} +
C_2(s, {\tilde y}_{2l-1}(s)) \, (s - {\tilde y}_{2l-1}(s))^{\beta_2(s)} \\ \nonumber
&= V((s - {\tilde y}_{2l-1}(s))-, s, {\tilde y}_{2l-1}(s)+; b(s, {\tilde y}_{2l-1}(s)+)) 
\end{align}
for $s > 0$ and $l = 1, \ldots, {\tilde l}'-1$, where the right-hand sides are given by
(\ref{V34p})-(\ref{Ci34p}) with $b(s, {\tilde y}_{2l-1}(s)+) = b(s, {\tilde y}_{2l}(s)) = s$.
Moreover, in the region ${\tilde R}_{2{\tilde l}'}$, the condition of (\ref{4condcall3}),
for $l = {\tilde l}'$, changes its form to $C_2(\eps, 0) \to 0$ as $\eps \downarrow 0$,
since otherwise $V(\eps, \eps, 0) \to \pm \infty$ as $\eps \downarrow 0$, that must be
excluded by virtue of the obvious fact that the value function in (\ref{V5b}) is bounded
at zero, while the condition of (\ref{4condcall1}) holds for $l = {\tilde l}'$ as well.

In addition, the process $(X, S, Y)$ can exit the region ${\tilde R}_{2l}$ in (\ref{tilR3_2lp})
passing to the stopping region $D_*$ from (\ref{D4}) only through the point
$({\bar s}(y), {\bar s}(y), y)$, by hitting the plane $d_1$, so that increasing its second
component $S$ until it reaches the value ${\bar s}(y) = \inf \{ q > s \, | \, b(q, y) \leq q \}$.
Then, the candidate value function should be continuous at the point $({\bar s}(y), {\bar s}(y), y)$,
that is expressed by the equality
\begin{align}
\label{4condcall2}
&C_1({\bar s}(y)-, y) \, ({\bar s}(y)-)^{\beta_1({\bar s}(y)-)}
+ C_2({\bar s}(y)-, y) \, ({\bar s}(y)-)^{\beta_2({\bar s}(y)-)} \\
\notag
&= V({\bar s}(y), {\bar s}(y), y; \, b({\bar s}(y), y)) \equiv {\bar s}(y) - K
\end{align}
for each ${\tilde y}_{2l}(s) < y \leq {\tilde y}_{2l-1}(s)$, $l = 1, \ldots, {\tilde l}'-1$.
However, in the region ${\tilde R}_{2{\tilde l}'}$, we have ${\bar s}(y) =
\infty$, since for the points $(x,s,y) \in {\tilde R}_{2{\tilde l}'}$ satisfying
$s-y\geq K$, we have $b(s,y)>s$ which will hold irrespective of how large $s$
becomes.

Thus, the condition of (\ref{4condcall2}) changes its form to
$C_1(\infty, y) = 0$, since otherwise $V(x, \infty, y) \to \pm \infty$ as $x \uparrow
\infty$, that must be excluded by virtue of the obvious fact that the value function in
(\ref{V5b}) is bounded at infinity. We can therefore conclude that the candidate value
function admits the representation
\begin{align}
\label{4Vint31}
&V(x, s, y; {\bar s}(y), {\tilde y}_{2l-1}(s), {\tilde y}_{2l}(s)) \\
\notag
&= C_1(s, y; {\bar s}(y), {\tilde y}_{2l-1}(s), {\tilde y}_{2l}(s)) \, x^{\beta_1(s)}
+ C_2(s, y; {\bar s}(y), {\tilde y}_{2l-1}(s), {\tilde y}_{2l}(s)) \, x^{\beta_2(s)}
\end{align}
in the regions ${\tilde R}_{2l}$ given by (\ref{tilR3_2lp}), where
$C_i(s, y; {\bar s}(y), {\tilde y}_{2l-1}(s), {\tilde y}_{2l}(s))$, $i = 1, 2$, provide a
solution of the two-dimensional coupled system of first-order linear partial differential equations in
(\ref{B31c})-(\ref{B31d}) with the boundary conditions of (\ref{4condcall1})-(\ref{4condcall2}),
for $l = 1, \ldots, {\tilde l}'$.

In order to argue the existence and uniqueness of solutions of the
boundary value problem formulated above, let us use the classical results of
the general theory of linear systems of first-order partial differential equations
(see, e.g. \cite[Chapter~I]{DiB} or \cite[Chapter~VII]{Han}).
For this, we first observe that the system of first-order linear partial
differential equations in (\ref{B31c})-(\ref{B31d}) does not admit characteristic curves,
so that the considered system is of elliptic type. In this respect, we can consider
an invertible analytic function $y = \Upsilon(s)$ such that the corresponding non-characteristic
curve coincides with the boundary curve
$y_{2l-1}(s)$ or $y_{2l}(s)$ defined above
for some $l = 1, \ldots, l'$.
Then, we may apply an appropriate affine transformation, which takes the point
$(s, y_{2l-1}(s))$ or $(s, y_{2l}(s))$
into the origin, and introduce the change of coordinates from $(s, y)$ to
$(s, z)$ with $z = y - \Upsilon(s)$, in order to reduce the system of
(\ref{B31c})-(\ref{B31d}) to the normal form.
On the other hand, we can consider an invertible analytic
function $s = \Gamma(y)$ such that the corresponding non-characteristic curve coincides with
the boundary curve
$s(y)$ defined above.
In that case, we may apply an appropriate affine transformation, which
takes the point $(s(y), y)$
into the origin, and introduce the change of
coordinates from $(s, y)$ to $(q, y)$ with $q = s - \Gamma(y)$, in order
to reduce the system of (\ref{B31c})-(\ref{B31d}) to the normal form.
In both cases, taking into account the assumption of continuity of the partial
derivatives of $\delta(s, y)$ and $\sigma(s, y)$ on $[0, \infty]^2$,
we can conclude by means of a version of the Cauchy-Kowalewski theorem
from \cite[Chapter~I, Theorem~5.1]{DiB} or \cite[Theorem~7.2.9]{Han}
(also in connection with Holmgren's uniqueness theorem)
that there exists a (locally) unique solution of the system
(\ref{B31c})-(\ref{B31d}), satisfying the boundary conditions of
(\ref{4condcall1})-(\ref{4condcall2}).
The obtained solution can admit an analytic continuation into the appropriate
parts of the state space $E^3$
(see, e.g. \cite{Z}, \cite{U}, \cite{Malek}, \cite{Miyake} and the references therein).

\def\skip{
Here and after, the existence and uniqueness of solutions to such systems
of this special kind follow from the classical results of the
existence and uniqueness of solutions to the appropriate boundary value
problems for two-dimensional linear systems of
More precisely, we consider non-characteristic analytic curves about
the boundary points in the conditions (\ref{4condcall1})-(\ref{4condcall2})
above, or (\ref{condcall1})-(\ref{condcall2}) below. Then, we introduce the
change of the coordinates from $(s, y)$ to $(s, z)$ with $z = y - \Upsilon(s)$
or from $(s, y)$ to $(q, y)$ with $q = s - \Gamma(y)$, for some analytic
functions $\Upsilon(s)$ and $\Gamma(y)$, to reduce the system of
(\ref{B31c})-(\ref{B31d}) to the normal form and such that
the functions $\Upsilon(s)$ and $\Gamma(y)$ are chosen with the
aim to coincide with the boundary functions $y_{2l-1}(s)$ or
$y_{2l}(s)$ for some $l \in \NN$ and $s(y)$ defined above.
Hence, taking into account the assumption of continuity of the partial
derivatives of $\delta(s, y)$ and $\sigma(s, y)$ on $[0, \infty]^2$,
and applying (locally) the appropriate affine transformations of the coordinates,
we can flatten the curves $y = \Upsilon(s)$ and $s = \Gamma(y)$ about the origin
}

Note that such coupled systems of first-order linear partial differential
equations have recently arisen in \cite{Pe5b} under
the study of certain other optimal stopping problems for the running extremal processes.
However, we may observe that the system in (\ref{B31c})-(\ref{B31d}) above turns out
to be essentially more complicated than the corresponding system in \cite[Equations~(3.42)-(3.43)]{Pe5b},
because the latter system can be decoupled.
The difficulty for the former system also arises from the specific form of the boundary conditions of (\ref{4condcall1})-(\ref{4condcall2})
formulated above or (\ref{condcall1})-(\ref{condcall2}) below.
The complicated structure of these conditions can be explained by the fact that
the running maximum process $S$ plays a crucial role in the definition of the
running maximum drawdown process $Y$ in (\ref{S4}), but not in the definition
of the running minimum process $I$, which is the counterpart coordinate process
contained in
\cite{Pe5b}.


\vspace{3pt}

{\bf (ii) The case of floating strike.}
Let us now consider the payoff function $G(x,s,y) = (Kx-s+y)^+$ in (\ref{V5b}).
Then, solving the system of equations in (\ref{B31b})-(\ref{B31bb}), we obtain
that the function in (\ref{V0}) admits the representation
\begin{equation}
\label{V35a}
V(x, s, y; b(s, y)) = C_1(s, y; b(s, y)) \, x^{\beta_1(s)} +
C_2(s, y; b(s, y)) \, x^{\beta_2(s)}
\end{equation}
for $0 < s-y \le x < b(s, y) \leq s$, with
\begin{equation}
\label{Ci35a}
C_i(s, y; b(s, y)) =
\frac{(\beta_{3-i}(s) - 1) K b(s, y) - \beta_{3-i}(s) (s - y)}
{(\beta_{3-i}(s) - \beta_i(s)) b(s, y)^{\beta_i(s)}}
\end{equation}
for all $0 < y < s$ and $i = 1, 2$.
Hence, assuming that the boundary function $b(s, y)$ is continuously
differentiable, we apply the condition of (\ref{B31d}) to the functions $C_i(s, y) =
C_i(s, y; b(s, y))$, $i = 1, 2$, in (\ref{Ci35a}) to obtain that the boundary
solves the following first-order nonlinear ordinary differential equation
\begin{align}
\label{g'35a}
\partial_y b(s, y) &= \sum_{i=1}^{2} \frac{\beta_{3-i}(s) b(s, y)}
{(\beta_i(s) - 1)(\beta_{3-i}(s) - 1) K b(s, y) - \beta_i(s) \beta_{3-i}(s) (s-y)} \\
&\phantom{= \sum_{i=1}^{2} \;\:} \times \frac{((s - y)/b(s, y))^{\beta_{i}(s)}}
{((s - y)/b(s, y))^{\beta_{i}(s)} - ((s - y)/b(s, y))^{\beta_{3-i}(s)}} \notag
\end{align}
for $0 < y < s$.
Taking into account the condition in (\ref{WVW9}) for the value function in
(\ref{V35a})-(\ref{Ci35a}), we conclude that $b(s, y) \sim g_*(s) (s-y)$
should hold as $y \uparrow s$.
Here, the function $g_*(s)$ is the unique solution of the arithmetic equation
\begin{equation}
\label{g35a}
\sum_{i=1}^{2} (-1)^{i} \, \Big(\beta_i(s) (\beta_{3-i}(s) - 1)
- (\beta_i(s) - 1) (\beta_{3-i}(s) - 1) K g(s) \Big) \, g^{\beta_i(s)}(s) = 0
\end{equation}
so that any candidate solution of the differential equation in (\ref{g'35a})
should satisfy the condition
\begin{equation}
\label{startfloat}
\lim_{y \uparrow s} \frac{b(s, y)}{s-y} = g_*(s)
\end{equation}
for each $s > 0$.
(The proof of uniqueness of the solution of the equation in (\ref{g35a}) is given in the Appendix.)

\vspace{3pt}

\begin{picture}(160,108)
\put(20,15){\begin{picture}(120,92)

\put(0,0){\line(1,0){120}} \put(120,0){\line(0,1){92}}
\put(0,0){\line(0,1){92}} \put(0,92){\line(1,0){120}}   



\put(10,10){\vector(1,0){100}}
\put(10,10){\vector(0,1){80}}    

\put(81,9){\line(0,1){72}} 

\put(80,6){$s$}

\put(9,81){\line(1,0){72}} 

\put(10,81){\line(1,-1){71}} 

\put(6,80){${s}$}

\put(63,40){${\bf .}$}
\put(63,39.9){${\bf .}$}	\put(63.1,39.9){${\bf .}$}	\put(62.9,39.9){${\bf .}$}
\put(63,39.8){${\bf .}$}	\put(63.1,39.8){${\bf .}$} 	\put(62.9,39.8){${\bf .}$}
\put(63,39.7){${\bf .}$}	\put(63.1,39.7){${\bf .}$}	\put(62.9,39.7){${\bf .}$}
\put(63,39.6){${\bf .}$}

\put(61,36){$(x,s,y)$}

\put(70,40){\line(-1,0){19}}
\put(70,39.9){\line(-1,0){19}}
\put(70,40.1){\line(-1,0){19}}

\put(51,40){\line(-1,1){10}}
\put(51,39.9){\line(-1,1){10}}
\put(51,40.1){\line(-1,1){10}}
\put(51,39.8){\line(-1,1){10}}
\put(51,40.2){\line(-1,1){10}}

\put(75,43){\line(-1,0){27}}
\put(75,42.9){\line(-1,0){27}}
\put(75,43.1){\line(-1,0){27}}

\put(60,46.5){\line(-1,0){15.5}}
\put(60,46.4){\line(-1,0){15.5}}
\put(60,46.6){\line(-1,0){15.5}}

\put(72,50){\line(-1,0){31}}
\put(72,49.9){\line(-1,0){31}}
\put(72,50.1){\line(-1,0){31}}

\put(10,81){\line(2,-1){71}}		
\put(10,81.1){\line(2,-1){71}}
\put(10,81.2){\line(2,-1){71}}
\put(10,80.9){\line(2,-1){71}}
\put(10,80.8){\line(2,-1){71}}

\qbezier[71](10,45)(45.5,45)(81,45)      
\put(.5,43.4){${\tilde y}_{1}(s)$}

\put(29,23){\vector(1,1){19}} \put(19,20){$x = s-y$}

\put(95.5,41){\vector(-4,-1){14}} \put(97.5,41){$x = s$}

\put(92,75){\vector(-2,-1){34}} \put(86,77){$x = g_*(s) \, (s-y)$}

\put(6,88){$y$} \put(108,6){$x$}
\end{picture}}
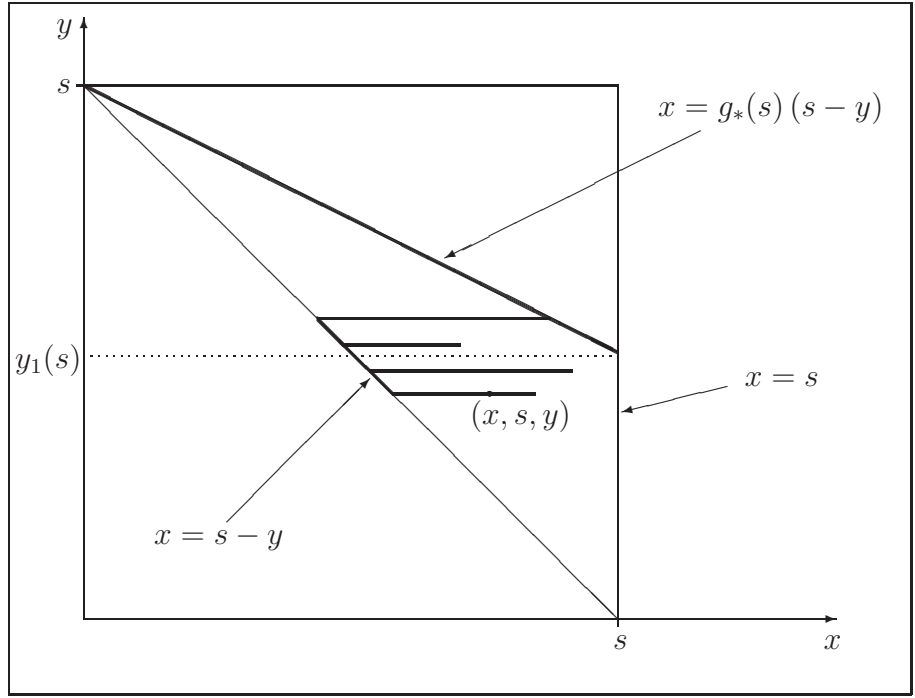
\put(25,10){\small{{\bf Figure 2.} A computer drawing of
$b_*(s, y)=g_*(s) (s-y)$ for}} \\
\put(25,5) {\small{some $s>0$ fixed, in the case when $\gamma_i(s,y) =
\beta_i(s)$, for $i=1,2$.}}
\end{picture}

For any $s > 0$ fixed, we observe that the ordinary differential equation in (\ref{g'35a})
is equivalent to a one with separable variables and admits the explicit solution $b_*(s, y) = g_*(s)(s - y)$
which satisfies the starting condition of (\ref{startfloat}) and
stays strictly above the plane $d_2$ and the surface $\{(x,s,y) \in E^3 \, | \, x = (s-y)/K \vee {\underline b}(s, y) \}$.
These assumptions for the boundary function $b_*(s, y)$ follow from the structure of
the continuation region $C_*$ in (\ref{C5}), which results to the expressions
of (\ref{b*float}) in Lemma \ref{lem}.
Then, we put ${\tilde y}_0(s) = s$, ${\tilde y}_1(s) = (s (g_*(s) - 1) / g_*(s))-$
and ${\tilde y}_2(s) = 0$, and observe that the boundary $b_*(s, y)$ exits the
region $E^3$ from the side of $d_1$ at the point $(s, s, {\tilde y}_1(s))$ and never returns back.
Hence, the candidate value function admits the expression in
(\ref{V35a})-(\ref{Ci35a}) in the region ${\tilde R}_{1}$ in
(\ref{tilR3_2l-1p}) and the boundary $b_*(s, y) = g_*(s) (s-y)$ provides
the explicit solution of the equation in (\ref{g'35a}), satisfying the 
condition of (\ref{startfloat}) and such that (\ref{b*float}) holds
(see Figure 2 above).

On the other hand, the candidate value function takes the form of (\ref{V0})
with $C_{i}(s, y)$, $i = 1, 2$, solving the linear system of first-order
partial differential equations in (\ref{B31c})-(\ref{B31d}),
in the region ${\tilde R}_{2}$ in (\ref{tilR3_2lp}), which belongs to the
continuation region $C_*$ in (\ref{C5}).
We can therefore conclude by means of the arguments presented in part (i)
above that the candidate value function admits the
representation (\ref{4Vint31}) for $l=1$ in the region ${\tilde R}_{2}$ given by
(\ref{tilR3_2lp}),
where $C_i(s, y; {\bar s}(y), {\tilde y}_{1}(s), {\tilde y}_{2}(s))$,
$i = 1, 2$, provide a unique solution of the two-dimensional system of first-order
linear partial differential equations in (\ref{B31c})-(\ref{B31d}) with the
boundary condition of (\ref{4condcall1}), where the right-hand side is
given by (\ref{V35a})-(\ref{Ci35a}) with $b(s, {\tilde y}_1(s)) = s$,
as well as the boundary conditions $C_2(\eps, 0) \to 0$ as $\eps \downarrow 0$
and $C_1(\infty, y) = 0$.

\vspace{6pt}


{\bf 3.3. The solution to the problem in the general setting.}
We now continue with the general form of the coefficients $\delta(s, y)$ and
$\sigma(s, y)$ in (\ref{dX4}), and thus, of the functions
$\gamma_i(s, y)$, $i = 1, 2$, from (\ref{gamma12}).

\vspace{3pt}

{\bf (i) The case of fixed strike.}
Let us now consider the payoff $G(x, s, y) = (K - s + y)^+$ in (\ref{V5b}).
In this case, solving the system of equations in (\ref{B31b})-(\ref{B31bb}),
we obtain that the function in (\ref{V0}) admits the representation
\begin{equation}
\label{V34}
V(x, s, y; b(s, y)) = C_1(s, y; b(s, y)) \, x^{\gamma_1(s, y)} +
C_2(s, y; b(s, y)) \, x^{\gamma_2(s, y)}
\end{equation}
for $0 < s-y \le x < b(s, y) \leq s$, with
\begin{equation}
\label{Ci34}
C_i(s, y; b(s, y)) = \frac{\gamma_{3-i}(s, y) ( K - s + y) }
{(\gamma_{3-i}(s, y) - \gamma_i(s, y)) b(s, y)^{\gamma_i(s, y)}}
\end{equation}
for all $0 < y < s$ and $i = 1, 2$.
Hence, assuming that the boundary function $b(s, y)$ is continuously differentiable,
we apply the condition of (\ref{B31d}) to the functions $C_i(s, y) =
C_i(s, y; b(s, y))$, $i = 1, 2$, in (\ref{Ci34}) to obtain that the boundary
solves the first-order nonlinear ordinary differential equation
\begin{align}
\label{g'34}
\partial_y b(s, y) &= \sum_{i=1}^{2} \frac{b(s, y)}{\gamma_i(s,y)}
\bigg( \frac{((s-y)/b(s, y))^{\gamma_i(s, y)}}{((s-y)/b(s, y))^{\gamma_{i}(s,y)}
- ((s-y)/ b(s, y))^{\gamma_{3-i}(s,y)}} \\
&\phantom{=\sum_{i=1}^{2} \;\:}
\times \bigg( \frac{1}{K-s+y} + \partial_y \gamma_i(s,y) \, \ln {\frac{s-y}
{b(s,y)}} \bigg) + \frac{\partial_y \gamma_i(s,y)}
{\gamma_{3-i}(s,y) - \gamma_i(s,y)} \bigg) \nonumber
\end{align}
for $0 < y < s$, where the partial derivatives $\partial_y \gamma_{i}(s, y)$,
$i = 1, 2$, are given by (\ref{gammay}) with (\ref{psi}).
Since the functions $\delta(s, y)$ and $\sigma(s, y)$ are assumed
to be continuously differentiable and bounded, it follows that
the limits $\delta(s, s-)$ and $\sigma(s, s-)$ exist for each $s > 0$.
Then, the limits $\gamma_i(s, s-)$ can be identified with the functions
$\beta_i(s)$, $i = 1, 2$, from the solution of the problem in the
particular setting considered in the previous subsection.
Taking into account the condition of (\ref{WVW5}) for the value function in
(\ref{V5b}), we conclude after some straightforward calculations that
$b(s, y) \sim g_*(s) (s-y)$ should hold as $y \uparrow s$, where $g_*(s)$
solves the equation in (\ref{gfixp}) with $\beta_i(s) = \gamma_i(s, s-)$,
that eventually yields the unique solution $g_*(s) = 1$ for all $s > 0$.
Thus, any candidate solution of the differential equation in (\ref{g'34})
should satisfy the starting condition of (\ref{startfixp}).

For any $s > 0$ fixed, let us now consider 
a candidate solution $b(s, y)$ of the equation in (\ref{g'34}), satisfying the starting condition of (\ref{startfixp}),
given that this solution satisfies the expressions in (\ref{b*fix}).
Then,
we define a decreasing sequence $({\tilde y}_n(s))_{n \in \NN}$
as in part (i) of the previous subsection.
Therefore, the candidate value function admits the expression of (\ref{V34})-(\ref{Ci34})
in the regions ${\tilde R}_{2l-1}$ from (\ref{tilR3_2l-1p}),
for $l = 1, \ldots, {\tilde l}'$.

On the other hand, the candidate value function takes the form of
(\ref{V0}) with $C_{i}(s, y)$, $i = 1, 2$, solving the linear system of
first-order partial differential equations in (\ref{B31c})-(\ref{B31d})
in the regions ${\tilde R}_{2l}$ from (\ref{tilR3_2lp}), for $l = 1,\ldots,
{\tilde l}'$, which belong to the continuation region $C' \cup C''$ given in (\ref{C41}) and (\ref{C42}).
Following arguments similar to the ones from part (i) of the previous subsection,
we obtain that the value function satisfies the conditions
\begin{align}
\label{condcall1}
&C_1(s, {\tilde y}_{2l}(s)+) \, ((s - {\tilde y}_{2l}(s))-)^{\gamma_1(s, {\tilde y}_{2l}(s)+)}
+ C_2(s, {\tilde y}_{2l}(s)+) \, ((s - {\tilde y}_{2l}(s))-)^{\gamma_2(s, {\tilde y}_{2l}(s)+)} \\ \nonumber
&= V(s - {\tilde y}_{2l}(s), s, {\tilde y}_{2l}(s); b_*(s, {\tilde y}_{2l}(s))) \\
\label{condcall3}
&C_1(s, {\tilde y}_{2l-1}(s)) \, (s - {\tilde y}_{2l-1}(s))^{\gamma_1(s, {\tilde y}_{2l-1}(s))} +
C_2(s, {\tilde y}_{2l-1}(s)) \, (s - {\tilde y}_{2l-1}(s))^{\gamma_2(s, {\tilde y}_{2l-1}(s))} \\ \nonumber
&= V((s - {\tilde y}_{2l-1}(s))-, s, {\tilde y}_{2l-1}(s)+; b(s, {\tilde y}_{2l-1}(s)+)) 
\end{align}
for $s > 0$, where the right-hand sides are given by (\ref{V34})-(\ref{Ci34}) with
$b(s, {\tilde y}_{2l-1}(s)+) = b(s, {\tilde y}_{2l}(s)) = s$, and
\begin{align}
\label{condcall2}
&C_1({\bar s}(y)-, y) \, ({\bar s}(y)-)^{\gamma_1({\bar s}(y)-, y)}
+ C_2({\bar s}(y)-, y) \, ({\bar s}(y)-)^{\gamma_2({\bar s}(y)-, y)} \\
\notag
&= V({\bar s}(y), {\bar s}(y), y; b({\bar s}(y), y)) \equiv {\bar s}(y) - K
\end{align}
for each ${\tilde y}_{2l}(s) < y \leq {\tilde y}_{2l-1}(s)$ and $l = 1, \ldots, {\tilde l}'-1$.
Moreover, we similarly obtain the condition in (\ref{condcall1}) together with
$C_2(\eps, 0) \rightarrow 0$ as $\eps \downarrow 0$, and $C_1(\infty, y) = 0$,
instead of (\ref{condcall3}) and (\ref{condcall2}), respectively, for $l = {\tilde l}'$.
We can therefore conclude that the candidate value function admits the representation
\begin{align}
\label{Vint31}
&V(x, s, y; {\bar s}(y), {\tilde y}_{2l-1}(s), {\tilde y}_{2l}(s)) \\
\notag
&= C_1(s, y; {\bar s}(y), {\tilde y}_{2l-1}(s), {\tilde y}_{2l}(s)) \, x^{\gamma_1(s, y)}
+ C_2(s, y; {\bar s}(y), {\tilde y}_{2l-1}(s), {\tilde y}_{2l}(s)) \, x^{\gamma_2(s, y)}
\end{align}
in the regions ${\tilde R}_{2l}$ given by (\ref{tilR3_2lp}), where
$C_i(s, y; {\bar s}(y), {\tilde y}_{2l-1}(s), {\tilde y}_{2l}(s))$,
$i = 1, 2$, provide a unique solution of the two-dimensional coupled system of
first-order linear partial differential equations in (\ref{B31c})-(\ref{B31d})
with the boundary conditions of (\ref{condcall1})-(\ref{condcall2}),
for $l = 1, \ldots, {\tilde l}'$. The existence and uniqueness
of the solution of the latter system follows from the arguments
presented in part (i) of the previous subsection.


\vspace{3pt}


{\bf (ii) The case of floating strike.}
Let us now consider the payoff $G(x, s, y) = (K x - s + y)^+$ in (\ref{V5b}).
Then, solving the system of equations in (\ref{B31b})-(\ref{B31bb}),
we obtain that the function in (\ref{V0}) admits the representation
\begin{equation}
\label{V35}
V(x, s, y; b(s, y)) = C_1(s, y; b(s, y)) \, x^{\gamma_1(s, y)} +
C_2(s, y; b(s, y)) \, x^{\gamma_2(s, y)}
\end{equation}
for $0 < s-y \le x < b(s, y) \leq s$, with
\begin{equation}
\label{Ci35}
C_i(s, y; b(s, y)) = \frac{(\gamma_{3-i}(s, y) - 1) K b(s, y) -
\gamma_{3-i}(s, y) (s - y)}{(\gamma_{3-i}(s, y) - \gamma_i(s, y))
{b(s, y)^{\gamma_i(s, y)}}}
\end{equation}
for all $0 < y < s$ and $i = 1, 2$.
Hence, assuming that the boundary function $b(s, y)$ is continuously differentiable,
we apply the condition of (\ref{B31d}) to the functions $C_i(s, y) = C_i(s, y; b(s,
y))$, $i=1,2$, in (\ref{Ci35}) to obtain that the boundary solves the first-order
nonlinear ordinary differential equation
\begin{align}
\label{g'35}
\partial_y b(s, y) &= \sum_{i=1}^{2} \bigg( \frac{\gamma_{3-i}(s, y) b(s, y) }{(\gamma_i(s, y) - 1)
(\gamma_{3-i}(s, y) - 1) K b(s, y) - \gamma_i(s, y) \gamma_{3-i}(s, y) (s-y)}  \\
&\phantom{= \sum_{i=1}^{2}  \bigg( \;} \times \frac{((s - y)/b(s, y))^{\gamma_{i}(s, y)}}
{((s - y)/b(s, y))^{\gamma_{i}(s, y)} - ((s - y) /b(s, y))^{\gamma_{3-i}(s, y)}} \nonumber \\
&\phantom{= \sum_{i=1}^{2}  \bigg( \;}+
\frac{ b(s, y) ((\gamma_{3-i}(s, y) - 1) K b(s, y) -
\gamma_{3-i}(s, y) (s-y))}{(\gamma_1(s, y) - 1) (\gamma_2(s, y)
- 1) K b(s, y) - \gamma_1(s, y) \gamma_2(s, y) (s-y)} \, \partial_y \gamma_{i}(s, y)
\nonumber \\
&\phantom{= \sum_{i=1}^{2}  \bigg( \;} \times \bigg( \frac{1}{\gamma_{3-i}(s,y) - \gamma_i(s,y)}
+ \frac{((s - y)/b(s, y))^{\gamma_{i}(s,y)} \ln{((s - y)/b(s,y))}}
{((s - y)/b(s, y))^{\gamma_{i}(s, y)}-((s - y)/b(s, y))^{\gamma_{3-i}(s, y)}} \bigg)
\bigg) \nonumber
\end{align}
for $0 < y < s$, where the partial derivatives $\partial_y \gamma_{i}(s, y)$,
$i = 1, 2$, are given by (\ref{gammay}) with (\ref{psi}).
Recall the assumption that the functions $\delta(s, y)$ and $\sigma(s, y)$ are
continuously differentiable and bounded,
so that the limits $\gamma_i(s, s-)$ can be identified with the functions
$\beta_i(s)$, $i = 1, 2$, from the solution of the problem in the
particular setting considered in the previous subsection.
Therefore, the function in (\ref{V35})-(\ref{Ci35}) should satisfy the property
$V(x, s, y; b(s, y)) \rightarrow V(x, s, s-\eps; b(s, s-\eps))$ as $y \uparrow s - \eps$,
for each $s-y \leq x < b(s, y)$ and any sufficiently small $\eps > 0$,
where $V(x, s, s-\eps; b(s, s-\eps))$ is given by the equation in
(\ref{V35a})-(\ref{Ci35a}), for $\eps < x < b(s, s-\eps)$, with $b(s, s-\eps)$
being a solution of the differential equation in (\ref{g'35a}).
Thus, letting $\eps \downarrow 0$, we see that any candidate solution of the
differential equation in (\ref{g'35}) should satisfy the starting condition of (\ref{startfloat}).

For any $s > 0$ fixed, let us now consider a candidate solution $b(s, y)$ of the equation
in (\ref{g'35}), satisfying the starting condition of (\ref{startfloat}),
given that this solution satisfies the expressions in (\ref{b*float}).
Then,
we define a decreasing sequence $({\tilde y}_n(s))_{n \in \NN}$ as in part (i) of
the previous subsection.
Therefore, the candidate value function admits the expression in
(\ref{V35})-(\ref{Ci35}) in the regions ${\tilde R}_{2l-1}$ defined in
(\ref{tilR3_2l-1p})
for $l = 1, \ldots, {\tilde l}'$.

On the other hand, the candidate value function takes the form of (\ref{V0})
with $C_{i}(s, y)$, $i = 1, 2$, solving the linear system of first-order
partial differential equations in (\ref{B31c})-(\ref{B31d}) in the
regions ${\tilde R}_{2l}$ defined in (\ref{tilR3_2lp}),
for $l = 1, \ldots, {\tilde l}'$, which belong to the continuation region $C_*$ in (\ref{C5}).
Using arguments similar to the ones of part (i) of this subsection, we can obtain the
same conditions as in (\ref{condcall1})-(\ref{condcall2}) above, where the
right-hand sides are given by (\ref{V35})-(\ref{Ci35}) with
$b(s, {\tilde y}_{2l-1}(s)+) = b(s, {\tilde y}_{2l}(s)) = s$.
By means of the arguments from part (i) of the previous subsection,
we can therefore conclude that the candidate value function admits the representation
in (\ref{Vint31}) in the regions ${\tilde R}_{2l}$ given by (\ref{tilR3_2lp}),
where $C_i(s, y; {\bar s}(y), {\tilde y}_{2l-1}(s), {\tilde y}_{2l}(s))$,
$i = 1, 2$, provide a unique solution of the two-dimensional coupled system of
first-order linear partial differential equations in (\ref{B31c})-(\ref{B31d})
satisfying the boundary conditions of (\ref{condcall1})-(\ref{condcall2}),
for $l = 1, \ldots, {\tilde l}'$.

     \section{\dot Main result and proof}

     In this section, we formulate and prove the main result of the paper, using the facts proved above.
     The proof of this assertion is based on a development of the maximality principle established
     in \cite{Pmax} and its extension to an optimal stopping problem for a three-dimensional Markov
     process $(X, S, Y)$ from (\ref{dX4})-(\ref{S4}) (see also \cite{Pe5b} for another three-dimensional problem).

\begin{theorem} \label{thmSY}
In the perpetual American fixed-strike or floating-strike option on
the maximum of market depth with payoff $G(x,s,y) = (K-s+y)^+$ or
$G(x,s,y) = (Kx-s+y)^+$, the value function of the optimal stopping
problem of (\ref{V5b}) for the process $(X, S, Y)$ from (\ref{dX4})-(\ref{S4})
has the expression
\begin{equation}
\label{V*42}
V_*(x, s, y)=
\begin{cases}
V(x, s, y; b_*(s, y)), & \text{if} \quad s - y \leq x < b_*(s, y) \leq s \\
V(x, s, y; {\bar s}(y), {\tilde y}_{2l-1}(s), {\tilde y}_{2l}(s)), & \text{if} \quad s-y \leq x \le s < b_*(s, y) \\
G(x, s, y), & \text{if} \quad (s-y) \vee b_*(s, y) \leq x \le s
\end{cases}
\end{equation}
and the optimal stopping time is given by (\ref{tau*}), where the functions
$V(x, s, y; b_*(s, y))$ and $V(x, s, y; {\bar s}(y), {\tilde y}_{2l-1}(s), {\tilde y}_{2l}(s))$
as well as the boundary function $b_*(s, y)$ are specified as follows:

(i) if $G(x, s, y) = (K - s + y)^+$ then
the function $V(x, s, y; b_*(s, y))$ is given by (\ref{V34})-(\ref{Ci34}) and the
boundary $b_*(s, y)$ provides the minimal solution of the equation in (\ref{g'34})
satisfying the starting condition of (\ref{startfixp}) and such that (\ref{b*fix}) holds for
$(x, s, y) \in {\tilde R}_{2l-1}$ defined in (\ref{tilR3_2l-1p}), and $V(x, s, y;
{\bar s}(y), {\tilde y}_{2l-1}(s), {\tilde y}_{2l}(s))$
is given by (\ref{Vint31}), whenever $(x, s, y) \in {\tilde R}_{2l}$ defined in
(\ref{tilR3_2lp}), with $C_i(s, y)$, $i = 1, 2$, solving the coupled system of equations in
(\ref{B31c})-(\ref{B31d}) and satisfying the conditions of (\ref{condcall1})-(\ref{condcall2}),
$l = 1, \ldots, {\tilde l}'$, where (\ref{condcall3})-(\ref{condcall2}) change
their form to $C_2(\eps, 0) \rightarrow 0$ as $\eps \downarrow 0$, and $C_1(\infty,
y) = 0$, for the case $l = {\tilde l}'$;

(ii) if $G(x, s, y) = (K x - s + y)^+$ then
the function $V(x, s, y; b_*(s, y))$ is given by (\ref{V35})-(\ref{Ci35}) and the
boundary $b_*(s, y)$ provides the minimal solution of the equation in (\ref{g'35})
satisfying the starting condition of (\ref{startfloat}) and such that (\ref{b*float}) holds
for $(x, s, y) \in {\tilde R}_{2l-1}$ defined in (\ref{tilR3_2l-1p}), and $V(x, s,
y; {\bar s}(y), {\tilde y}_{2l-1}(s), {\tilde y}_{2l}(s))$ is given by (\ref{Vint31}),
whenever $(x, s, y) \in {\tilde R}_{2l}$ defined in (\ref{tilR3_2lp}), with $C_i(s, y)$, $i = 1, 2$,
solving the coupled system of equations in (\ref{B31c})-(\ref{B31d}) and satisfying the conditions of
(\ref{condcall1})-(\ref{condcall2}), $l = 1, \ldots, {\tilde l}'$, where (\ref{condcall3})-(\ref{condcall2})
change their form to $C_2(\eps, 0) \rightarrow 0$ as $\eps \downarrow 0$, and $C_1(\infty,
y) = 0$, for the case $l = {\tilde l}'$;

(iii) if $G(x, s, y) = (K x - s + y)^+$ as well as $\delta(s,y) = \delta(s)$ and
$\sigma(s,y) = \sigma(s)$ 
then the function
$V(x, s, y; b_*(s, y))$ is given by (\ref{V35a})-(\ref{Ci35a}) and the boundary
is given by $b_*(s, y) = g_*(s) (s-y)$
as the explicit solution of the equation in (\ref{g'35a}) satisfying the starting
condition of (\ref{startfloat}) and such that (\ref{b*float}) holds,
where $g_*(s)$ provides the unique solution of (\ref{g35a})
for $(x, s, y) \in {\tilde R}_{1}$ defined in (\ref{tilR3_2l-1p}), and
$V(x, s, y; {\bar s}(y), {\tilde y}_1(s), {\tilde y}_2(s))$
is given by (\ref{4Vint31}), whenever $(x, s, y) \in {\tilde R}_{2}$
defined in (\ref{tilR3_2lp}), with $C_i(s, y)$, $i = 1, 2$, solving the coupled
system of equations in (\ref{B31c})-(\ref{B31d}) and satisfying the conditions of
(\ref{4condcall1}), $C_2(\eps, 0) \rightarrow 0$ as $\eps \downarrow 0$, and
$C_1(\infty, y) = 0$, for $l = l' = 1$.
\end{theorem}

Since all the assertions formulated above are proved using similar arguments, we only
give a proof for the general optimal stopping problem related to the perpetual
American fixed-strike option on the maximum of market depth in part (i) of Theorem \ref{thmSY}.

      \vspace{6pt}

{\bf Proof of part (i) of Theorem 4.1.}
In order to verify the assertion stated above, it remains to show that the function
defined in (\ref{V*42}) coincides with the value function in (\ref{V5b}) with payoff
$(K - s + y)^+$ and that the stopping time $\tau_*$ in (\ref{tau*}) is optimal with
the boundary $b_*(s, y)$ specified above.
For this, let $b(s, y)$ be any solution of (\ref{g'34})
with the starting condition in (\ref{startfixp}) and satisfying (\ref{b*fix}).
Let us also denote by $V_{b}(x, s, y)$ the right-hand side of the expression
in (\ref{V*42}) associated with this $b(s, y)$.
It then follows using straightforward calculations and the assumptions presented
above that the function $V_{b}(x, s, y)$ solves the system of (\ref{LV})-(\ref{VD}),
while the normal-reflection and smooth-fit conditions are satisfied in (\ref{NRY})-(\ref{NRboth}).
Hence, taking into account the fact that the function $V_{b}(x, s, y)$ is $C^{2,1,1}$
and the boundary $b(s, y)$ is assumed to be continuously differentiable for all $0 <
y < s$, by applying the change-of-variable formula from \cite[Theorem~3.1]{Pe1a} to
$e^{- r t} V_{b}(X_t, S_t, Y_t)$, we obtain
\begin{align}
\label{4rho4c}
&e^{- r t} \, V_{b}(X_t, S_t, Y_t) = V_{b}(x, s, y) + M_t \\
\notag
&+ \int_0^t e^{- r u} \, (\LL V_{b} - r V_{b}) (X_u, S_u, Y_u) \, I(X_u \not= S_{u} -
Y_u, X_u \not= b(S_u, Y_u), X_u \not= S_u) \, du \\
\notag
&+ \int_0^t e^{- r u} \, \partial_s V_{b}(X_u, S_u, Y_u) \, I(X_u = S_u) \, dS_u
+ \int_0^t e^{- r u} \, \partial_y V_{b}(X_u, S_u, Y_u) \, I(X_u = S_u - Y_u) \, dY_u
\end{align}
where 
the process $M = (M_t)_{t \ge 0}$ given by
\begin{equation}
\label{4N5}
M_t = \int_0^t e^{- r u} \, \partial_x V_{b}(X_u, S_u, Y_u) \, I(X_u \neq S_{u} -
Y_u, X_u \neq S_{u})
\, \sigma(S_u, Y_u) \, X_{u} \, dB_u
\end{equation}
is a square integrable martingale under $P_{x, s, y}$.
 Note that, since the time spent by the process $X$ at the boundary
surface $\{ (x, s, y) \in E^3 \, | \, x = b(s, y) \}$
as well as at the planes $d_1 = \{(x, s, y) \in \RR^3 \, | \, 0 < x = s \}$ 
and $d_2 = \{(x, s, y) \in \RR^3 \, | \, 0 < x = s-y \}$ 
is of Lebesgue measure zero, the indicators in the second line of the formula
(\ref{4rho4c}) as well as in the formula (\ref{4N5}) can be ignored.
Moreover, since the process $S$ increases only on the plane $d_1$
and the process $Y$ increases only on the plane $d_2$,
the indicators in the third line of (\ref{4rho4c}) can be set equal to one.

By using straightforward calculations and the arguments from the previous section, it
is verified that $(\LL V_{b} - r V_{b})(x, s, y) \le 0$ for all $(x, s, y) \in E^3$
such that $x \neq b(s, y)$, $x \neq s - y$, and $x \neq s$.
Moreover, it is shown by means of standard arguments that the properties in (\ref{VC})-(\ref{LVD})
also hold, which together with (\ref{CF})-(\ref{VD}) imply that the inequality
$V_{b}(x, s, y) \ge (K - s + y)^+$ is satisfied for all $(x, s, y) \in E^3$.
It therefore follows from the expression (\ref{4rho4c}) that the inequalities
\begin{equation}
\label{4rho4e}
e^{- r \tau} \, 
(K - S_{\tau} + Y_{\tau})^+
\le e^{- r \tau} \, V_{b}(X_{\tau}, S_{\tau}, Y_{\tau})
\le V_{b}(x, s, y) + M_{\tau}
\end{equation}
hold for any finite stopping time $\tau$ with respect to the natural filtration of $X$.

Taking the expectation with respect to $P_{x, s, y}$ in (\ref{4rho4e}), by means of
the optional sampling theorem (see, e.g.
\cite[Chapter~I, Theorem~3.22]{KS}), we get
\begin{align}
\label{4rho4e2}
E_{x, s, y} \big[ e^{- r (\tau \wedge t)} \, 
(K - S_{\tau \wedge t} + Y_{\tau \wedge t})^+ 
\big] &\le
E_{x, s, y} \big[ e^{- r (\tau \wedge t)} \, V_{b}(X_{\tau \wedge t}, S_{\tau \wedge
t}, Y_{\tau \wedge t}) \big] \\
\notag
&\le V_{b}(x, s, y) + E_{x, s, y} \big[ M_{\tau \wedge t} \big] = V_{b}(x, s, y)
\end{align}
for all $(x, s, y) \in E^3$.
Hence, letting $t$ go to infinity and using Fatou's lemma, we obtain that the inequalities
\begin{equation}
\label{4rho4e3}
E_{x, s, y} \big[ e^{- r \tau} \, 
(K - S_{\tau} + Y_{\tau})^+
\big] \le E_{x, s, y} \big[ e^{- r \tau} \, V_{b}(X_{\tau}, S_{\tau}, Y_{\tau}) \big]
\le V_{b}(x, s, y)
\end{equation}
are satisfied for any finite stopping time $\tau$ and all $(x, s, y) \in E^3$.
Taking first the supremum over all stopping times $\tau$ and then the infimum over
all $b$, we conclude that
\begin{equation}
\label{4rho4e32}
E_{x, s, y} \big[ e^{- r \tau_*} \, (K - S_{\tau_*} + Y_{\tau_*})^+ \big]
\le \inf_{b} V_{b}(x, s, y) = V_{b_*}(x, s, y)
\end{equation}
where $b_*(s, y)$ is the minimal solution of (\ref{g'34}) with the
starting condition in (\ref{startfixp}) and satisfying (\ref{b*fix}).
Using the fact that the function $V_{b}(x, s, y)$ is increasing in $b$,
satisfying $b(s, y) > s-y$ for $s-y < K$ and $b(s, y) > s$ for $s-y \geq K$
under any $0 < y < s$ fixed, we see that the infimum in (\ref{4rho4e32}) is
attained over any sequence of solutions $(b_n(s, y))_{n \in \NN}$ to (\ref{g'34})
with the starting condition in (\ref{startfixp}) and satisfying (\ref{b*fix}),
and such that $b_n(s, y) \downarrow b_*(s, y)$ as $n \rightarrow \infty$.
Since the inequalities in (\ref{4rho4e3}) hold also for $b_*(s, y)$, we see that
(\ref{4rho4e32}) holds for $b_*(s,y)$ and $(x,s,y) \in E^3$ as well.
Note that $V_{b}(x, s, y)$ in (\ref{4rho4e2}) is superharmonic
for the Markov process $(X,S,Y)$ on $E^3$.
Taking into account the fact that $V_{b}(x, s, y)$ is increasing in $b$ and that the
inequality $V_{b}(x, s, y) \geq (K-s+y)^+$ holds for all $(x,s,y) \in E^3$, we
observe that the selection of the minimal solution $b_*(s, y)$, which satisfies
(\ref{b*fix}) whenever such a choice exists, is equivalent to invoking the
superharmonic characterisation of the value function as the smallest superharmonic function
dominating the payoff function (see, e.g. \cite{Pmax} or \cite[Chapter~I, Section~2]{PSbook}).

In order to clarify the (local) existence and uniqueness of the solution to the
equation in (\ref{g'34}),
we recall the arguments of \cite[Subsection~3.5]{Pe5b} and denote the right-hand side
of (\ref{g'34}) by $\Psi(s, y, b(s, y))$. Since the function $\Psi(s, y, b)$ is (locally)
continuous in $y$ and (locally) Lipschitz in $b$, we may conclude from the general theory
of first-order nonlinear ordinary differential equations that the one in (\ref{g'34}) admits
a (locally) unique solution.

In order to construct the minimal solution $b_*(s, y)$ of the equation in (\ref{g'34}),
satisfying the conditions mentioned above, and prove the fact that it is optimal in $E^3$,
we provide an extension of the arguments from \cite[Theorem~3.1]{Pmax} to the $(X, S, Y)$-setting
(see also \cite[Subsection~3.5]{Pe5b} for another three-dimensional case).
For this, let us consider the sequence of stopping times $\tau_n$ defined as in (\ref{tau*})
with $b_n(s, y)$ instead of $b_*(s, y)$, where $b_n(s, y)$ is the solution of (\ref{g'34}) with
the starting condition of (\ref{startfixp}) and such that $b_n(s, y_n) = s - y_n < K$ holds
for some $y_n \downarrow c$ as $n \rightarrow \infty$, where $c > 0$ is such that
$b_\infty(s,y) > s-y$ for all $s-K < y \leq c$, and $b_\infty(s,y) > s$ holds for all
$y \leq s-K$, whenever such a sequence exists. Otherwise, we consider as $b_n(s, y)$
the solution of (\ref{g'34}) with the starting condition of (\ref{startfixp}) and such that
$b_n(s, s-K) = c_n$ holds for some $c_n \uparrow c$ as $n \rightarrow \infty$, where
$c > 0$ is such that $b_\infty(s,y) > s$ holds for all $y \leq s-K$.
It follows from the uniqueness of the solution to (\ref{g'34}) that no distinct
solutions intersect, so that the sequence $(b_n(s, y))_{n \in \NN}$ is increasing and
the limit $b_*(s, y) = \lim_{n \rightarrow \infty} b_n(s, y)$ exists.
By virtue of the fact that the function $V_{b_n}(x, s, y)$ from the expression in
(\ref{V*42}) associated with this $b_n(s, y)$, satisfies the system of
(\ref{LV})-(\ref{LVD}) with (\ref{NRY}) and taking into account the structure of
$\tau_n$ given by (\ref{tau*}) with $b_n(s, y)$ instead of $b_*(s, y)$, it follows from
the equivalent expression of (\ref{4rho4c}) that the equalities
\begin{align}
\label{4rho4g}
e^{- r (\tau_n \wedge t)} \, 
(K - S_{\tau_n \wedge t} + Y_{\tau_n \wedge t})^+
= e^{- r (\tau_n \wedge t)} \, V_{b_n}(X_{\tau_n \wedge t}, S_{\tau_n \wedge t}, Y_{\tau_n \wedge t})
= V_{b_n}(x, s, y) + M_{\tau_n \wedge t}
\end{align}
hold for all $(x, s, y) \in E^3$.
Observe that $\tau_n \uparrow \tau_*$ ($P_{x, s, y}$-a.s.),
the property
\begin{equation}
\label{ESS4}
E_{x, s, y} \Big[ \sup_{t \ge 0} e^{- \rho (\tau_* \wedge t)} \, Y_{\tau_* \wedge t} \Big]
\leq E_{x, s, y} \Big[ \sup_{t \ge 0} e^{- \rho (\tau_* \wedge t)} \, S_{\tau_* \wedge t} \Big]
= E_{x, s, y} \Big[ \sup_{t \ge 0} e^{- \rho (\tau_* \wedge t)} \, X_{\tau_* \wedge t} \Big] < \infty
\end{equation}
holds for all $(x, s, y) \in E$ and the variable $e^{- r \tau_*} (K - S_{\tau_*} + Y_{\tau_*})^+$
is bounded on the set $\{ \tau_* = \infty \}$.
Note also that, by using the asymptotic behavior of $b_*(s, y)$ when $y \uparrow s$ in
(\ref{startfixp}), the property $P_{x, s, y}(\tau_* < \infty) = 1$ holds for all
$(x, s, y) \in E^3$.
Hence, letting $t$ and $n$ go to infinity and using the conditions in (\ref{CF}) and
(\ref{NRY}), as well as the fact that $\tau_n \uparrow \tau_*$ ($P_{x, s, y}$-a.s.),
we can apply the Lebesgue dominated convergence theorem for (\ref{4rho4g}) to obtain the equality
\begin{equation}
\label{4rho4i}
E_{x, s, y} \big[ e^{- r \tau_*} \, 
(K - S_{\tau_*} + Y_{\tau_*})^+
\big] = V_{b_*}(x, s, y)
\end{equation}
for all $(x, s, y) \in E^3$, which together with (\ref{4rho4e32}) directly implies the desired assertion. $\square$


\section{Appendix}

In this section, we prove the existence of a unique solution $g_*(s)$ to the equation in
(\ref{g35a}). For this, we first rewrite the latter equation in the form
\begin{equation}
\label{I12gs}
F_1(g(s)) = F_2(g(s)) \;\;\; \text{with} \;\;\;
F_i(x) = \big((\beta_1(s) - 1) \, (\beta_2(s) - 1) \, K \, x -
\beta_i(s) \, (\beta_{3-i}(s) - 1)\big) \, x^{\beta_i(s)}
\end{equation}
for $i = 1, 2$ and some arbitrary function $g(s)$, for some $s>0$ fixed.
Then, the derivatives of the functions $F_i(x)$ from (\ref{I12gs})
take the form
\begin{align}
\label{Ii'}
F_i'(x) = (\beta_{3-i}(s) - 1) \, \big( (\beta^2_i(s) - 1) \, K \, x -
\beta^2_i(s) \big) \, x^{\beta_i(s) - 1}
\end{align}
for all $x>0$ and $i = 1, 2$, and some $s > 0$ fixed.
It therefore follows that
the function $F_1(x)$ is increasing on the interval
$(0, \beta^2_1(s)/((\beta^2_1(s) - 1) K))$ with
$F_1(0+) = 0$ and
\begin{equation} \label{I1u}
F_1 \bigg( \frac{\beta^2_1(s)}{(\beta^2_1(s) - 1) K}
\bigg) = - \frac{\beta_1(s) (\beta_2(s) - 1)}{\beta_1(s) + 1}
\bigg( \frac{\beta^2_1(s)}{(\beta^2_1(s) - 1) K} \bigg)^{\beta_1(s)} > 0
\end{equation}
and then decreasing on the interval $(\beta^2_1(s)/((\beta^2_1(s) - 1) K),
\infty)$ with $F_1(\infty) = - \infty$, for any $s > 0$ fixed.

Let us now specify the structure of the function $F_2(x)$ based on the value
of $\beta_2(s)$. For this, we first assume that $\beta_2(s) < -1$ holds, that is
equivalent to $2 r - \delta(s) - \sigma^2(s) > 0$ in the model in
(\ref{dX4})-(\ref{S4}), for some $s > 0$ fixed.
In this case, it follows that $F_2(x)$ is decreasing
on the interval $(0, \beta^2_2(s)/((\beta^2_2(s) - 1) K))$, with
$F_2(0+) = + \infty$ and
\begin{equation} \label{I2u}
F_2 \bigg( \frac{\beta^2_2(s)}{(\beta^2_2(s) - 1) K} \bigg) = - \frac{\beta_2(s)
(\beta_1(s) - 1)}{\beta_2(s) + 1}  \bigg( \frac{\beta^2_2(s)}{(\beta^2_2(s) - 1) K} \bigg)^{\beta_2(s)} < 0
\end{equation}
and increasing
on the interval $(\beta^2_2(s)/((\beta^2_2(s) - 1) K), \infty)$, with
$F_2(\infty) = 0$.
Thus, taking also into account the fact that $h_1(s) > h_2(s)$ holds as well,
where ${h_i}(s)$ is such that $F_i(h_i(s)) = 0$ for $i = 1, 2$, we conclude
that there exist exactly two solutions of the equation in (\ref{I12gs}).
Let us now assume that $-1 \leq \beta_2(s) < 0$ holds, that is equivalent to
$2 r - \delta(s) - \sigma^2(s) \leq 0$ in the model of (\ref{dX4})-(\ref{S4}).
In this case, the value $\beta^2_2(s) / (\beta^2_2(s) - 1)$ is negative, that
yields the fact that the function $F_2(x)$ is strictly decreasing in the interval
$(0, \infty)$, with $F_2(0+) = \infty$ and $F_2(\infty) = - \infty$.
Hence, the same arguments as above guarantee the existence of at least one solution
of the equation in (\ref{I12gs}). Moreover, by computing the second-order derivatives
of the functions $F_i(x)$, $i = 1, 2$, we get
\begin{equation}
\label{Ii''}
F_i''(x) = \beta_i(s) \, (\beta_1(s) - 1) \,
(\beta_2(s) - 1) \, \big( (\beta_i(s) + 1) K x - \beta_i(s) \big) \,
x^{\beta_i(s) - 2}
\end{equation}
for $s > 0$ fixed and $i = 1, 2$, we observe that $F_1''(x) < 0$
for all $x > h_1(s) > \beta^2_1(s)/((\beta^2_1(s) - 1) K) > \beta_1(s)
/((\beta_1(s) + 1)K) > 0$ and $F_2''(x) > 0$ for all
$x > h_2(s)$, under the assumption that $\beta_2(s) \ge -1$ holds.
It follows that the function $F_1(x)$ is concave, when it becomes
negative for $x \geq h_1(s) > h_2(s)$, while the function
$F_2(x)$ is convex, when it becomes negative for $x \geq h_2(s)$.
It therefore follows that there exist exactly two solutions of the equation
in (\ref{I12gs}).

Let us finally consider the two solutions $g_1(s)$ and $g_2(s)$ of the equation in
(\ref{I12gs}), which satisfy the inequalities $g_2(s) < h_2(s) < 1/K < h_1(s) < g_1(s)$,
for all possible choices of the parameters of the model.
The fact that $b_*(s,y)$ satisfies (\ref{b*float}) in Lemma 2.1, for all
$0<y<s$, implies that
\begin{equation}
\label{Ii'''}
g_*(s) > \frac{\rho}{\delta(s) K}
\equiv \frac{\beta_1(s) \beta_2(s)}{(\beta_1(s) - 1) (\beta_2(s) - 1) K}
> \frac{\beta_2(s)}{(\beta_2(s) - 1) K} \equiv h_2(s)
\end{equation}
holds for all $s > 0$.
We therefore reject $g_2(s)< h_2(s)$ and accept $g_1(s) > h_1(s) > \rho / (\delta(s) K)$.
Furthermore, taking into account the fact that $F_1(1) > F_2(1)$ and
$F_1(x) \leq F_2(x)$ holds for all $x \geq g_*(s)$, it follows that
the function $g_*(s) \equiv g_1(s)$ satisfying $F_1(g_*(s)) = F_2(g_*(s))$,
is greater than $1$, for any choice of fixed $s > 0$.

\vspace{12pt}

{\bf Acknowledgments.}
The authors are grateful to the Associate Referee and two anonymous Referees for their valuable suggestions which helped
to essentially improve the presentation of the paper. This research was supported by a Small Grant from the Suntory and Toyota International Centres for Economics and Related Disciplines (STICERD) at the London School of Economics and Political Science.


\begin{thebibliography}{MDS}






    \bibitem{AAP}
    {\sc Asmussen,~S.,~Avram,~F.} and {\sc Pistorius,~M.~(2003).}
    Russian and American put options under exponential phase-type
    L{\'e}vy models. {\it Stochastic Processes and their Applications} {\bf 109} (79--111).

    \bibitem{AKP}
    {\sc Avram,~F.,~Kyprianou,~A.~E.} and {\sc Pistorius,~M.~(2004).}
    Exit problems for spectrally negative L{\'e}vy processes and applications to (Canadized) Russian options.
    {\it Annals of Applied Probability} {\bf 14}(1) (215--238).

 \bibitem{BK3}
 {\sc Baurdoux,~E.~J.} and {\sc Kyprianou,~A.~E.~(2009).} The Shepp-Shiryaev
 stochastic game driven by a spectrally negative L{\'e}vy process. {\it
 Theory of Probability and its Applications} {\bf 53} (481--499).

\bibitem{BKP}
{\sc Baurdoux,~E.~J.,~Kyprianou,~A.~E.} and {\sc Pardo,~J.~C.~(2011).}
The Gapeev-K\"uhn stochastic game driven by a spectrally positive L\'evy process.
{\it Stochastic Processes and their Applications} {\bf 121} (1266--1289).







\bibitem{Det}
{\sc Detemple,~J.~(2006).}
{\it American-Style Derivatives: Valuation and Computation.}
Chapman and Hall/CRC, Boca Raton.

\bibitem{DiB}
{\sc DiBenedetto,~E.~(2010).} {\it Partial Differential Equations.} (Second Edition)
Birkh\"auser, Basel.

    \bibitem{DSS}
    {\sc Dubins,~L.,~Shepp,~L.~A.} and {\sc Shiryaev,~A.~N.~(1993).}
    Optimal stopping rules and maximal inequalities for Bessel processes.
    {\it Theory of Probability and its Applications} {\bf 38}(2) (226--261).



    \bibitem{Dyn}
    {\sc Dynkin,~E.~B.~(1963).}
    The optimum choice of the instant for stopping a Markov process.
    {\it Soviet Mathematical Doklady} {\bf 4} (627--629).

\bibitem{Forde}
{\sc Forde,~M.~(2011).}
A diffusion-type process with a given joint law for the terminal level and supremum at an independent exponential time.
{\it Stochastic Processes and their Applications} {\bf 121} (2802--2817).

\bibitem{FPZ}
{\sc Forde,~M.,~Pogudin,~A.} and {\sc Zhang,~H.~(2013).}
Hitting times, occupation times, tri-variate laws and the forward
Kolmogorov equation for a one-dimensional diffusion with memory.
{\it Advances in Applied Probability} {\bf 45} (860--875).

    \bibitem{Gap}
    {\sc Gapeev,~P.~V.~(2007).} Discounted optimal stopping for maxima
    of some jump-diffusion processes. {\it Journal of Applied Probability}
    {\bf 44} (713--731).

\bibitem{GR3}
{\sc Gapeev,~P.~V.} and {\sc Rodosthenous,~N.~(2014).}
Optimal stopping problems in diffusion-type models with running maxima and drawdowns.
{\it Journal of Applied Probability} {\bf 51} (799--817).

\bibitem{GHP} {\sc Glover,~K., Hulley,~H.} and {\sc Peskir,~G.~(2013).}
Three-dimensional Brownian motion and the golden ratio rule.
{\it Annals of Applied Probability} {\bf 23} (895--922).

    \bibitem{GP1}
    {\sc Graversen,~S.~E.} and {\sc Peskir,~G.~(1998).}
    Optimal stopping and maximal inequalities for geometric
    Brownian motion. {\it Journal of Applied Probability} {\bf 35}(4) (856--872).

    \bibitem{GP2}
    {\sc Graversen,~S.~E.} and {\sc Peskir,~G.~(1998).}
    Optimal stopping and maximal inequalities for linear diffusions.
    {\it Journal of Theoretical Probability} {\bf 11} (259--277).

    \bibitem{GuoShepp}
    {\sc Guo,~X.} and {\sc Shepp,~L.~A.~(2001).}
    Some optimal stopping problems with nontrivial boundaries
    for pricing exotic options. {\it Journal of Applied Probability}
    {\bf 38}(3) (647--658).

    \bibitem{GuoZer}
    {\sc Guo,~X.} and {\sc Zervos,~M.~(2010).} $\pi$ options.
    {\it Stochastic Processes and their Applications} {\bf 120}(7) (1033-–1059).


\bibitem{Han}
{\sc Han,~Q.~(2011).}
{\it A Basic Course in Partial Differential Equations.} Graduate Studies in Mathematics {\bf 120}.
American Mathematical Society, Providence.


    \bibitem{Jmax}
    {\sc Jacka,~S.~D.~(1991).}
    Optimal stopping and best constants for Doob-like inequalities I:
    The case $p=1$. {\it Annals of Probability} {\bf 19} (1798--1821).

\bibitem{KS} {\sc Karatzas,~I.} and {\sc Shreve,~S.~E.~(1991).}
\textit{Brownian Motion and Stochastic Calculus.} (Second Edition)
Springer, New York.

    \bibitem{LS}
    {\sc Liptser,~R.~S.} and {\sc Shiryaev,~A.~N.~(2001).}
    {\it Statistics of Random Processes I.}
    (Second Edition, First Edition 1977) Springer, Berlin.



    \bibitem{Malek}
    {\sc Malek,~S.~(2006).}
    On analytic continuation of solutions of linear partial differential equations.
{\it Journal of Dynamical and Control Systems} {\bf 12} (79--86).


    \bibitem{MP}
    {\sc Mijatovi{\'c},~A.} and {\sc Pistorius,~M.~R.~(2012).}
On the drawdown of completely asymmetric L{\'e}vy processes.
{\it Stochastic Processes and their Applications} {\bf 22} (3812--3836).

 \bibitem{Miyake}
 {\sc Miyake,~M.~(1981).}
Global and local Goursat problems in a class of holomorphic or partially holomorphic functions.
{\it Journal of Differential Equations} {\bf 39} (445--463).

    \bibitem{KyprOtt}
    {\sc Kyprianou,~A.~E.} and {\sc Ott,~C.~(2014).}
    A capped optimal stopping problem for the maximum process.
    {\it Acta Applicandae Mathematicae} {\bf 129} (147--174).

    \bibitem{Ott}
    {\sc Ott,~C.~(2013).}
    Optimal stopping problems for the maximum process with upper and lower caps.
    {\it Annals of Applied Probability} {\bf 23} (2327--2356).

    \bibitem{Ott_thesis}
    {\sc Ott,~C.~(2013).}
{\it Optimal Stopping Problems for the Maximum Process.}
PhD thesis, University of Bath.


    \bibitem{Jesper}
    {\sc Pedersen,~J.~L.~(2000).}
    Discounted optimal stopping problems for the maximum process.
    {\it Journal of Applied Probability} {\bf 37}(4) (972--983).

    \bibitem{Pmax}
    {\sc Peskir,~G.~(1998).}
    Optimal stopping of the maximum process: The maximality principle.
    {\it Annals of Probability} {\bf 26} (1614--1640).

    \bibitem{Pe1a}
    {\sc Peskir,~G.~(2007).}
    A change-of-variable formula with local time on surfaces.
    {\it S\'eminaire de Probababilit\'e, Lecture Notes in Mathematics} {\bf 1899} (69--96).

\bibitem{Pe5a}
{\sc Peskir,~G.~(2012).}
Optimal detection of a hidden target: The median rule.
{\it Stochastic Processes and their Applications} {\bf 122} (2249--2263).

\bibitem{Pe5b}
{\sc Peskir,~G.~(2014).}
Quickest detection of a hidden target and extremal surfaces.
{\it Annals of Applied Probability} {\bf 24} (2340--2370).

    \bibitem{PSbook}
    {\sc Peskir,~G.} and {\sc Shiryaev,~A.~N.~(2006).} {\it Optimal
    Stopping and Free-Boundary Problems.} Birkh{\"a}user, Basel.


\bibitem{PVH}
{\sc Pospisil,~L.,~Vecer,~J.} and {\sc Hadjiliadis,~O.~(2009).}
Formulas for stopped diffusion processes with stopping times based on
drawdowns and drawups. {\it Stochastic Processes and their Applications}
{\bf 119}(8) (2563--2578).




    \bibitem{SS1}
    {\sc Shepp,~L.~A.} and {\sc Shiryaev,~A.~N.~(1993).}
    The Russian option: reduced regret.
    {\it Annals of Applied Probability} {\bf 3}(3) (631--640).


    \bibitem{SS2}
    {\sc Shepp,~L.~A.} and {\sc Shiryaev,~A.~N.~(1994).}
    A new look at the pricing of Russian options.
    {\it Theory Probability and its Applications} {\bf 39}(1) (103--119).



    \bibitem{FM}
    {\sc Shiryaev,~A.~N.~(1999).}
    {\it Essentials of Stochastic Finance.}
    World Scientific, Singapore.

    \bibitem{U}
    {\sc Uchida,~M.~(1993).}
Continuation of analytic solutions of linear differential equations up to convex conical singularities.
{\it Bulletin de la Soci\'et\'e Math\'ematique de France} {\bf 121} (133--152).

    \bibitem{Z}
    {\sc Zerner,~M.~(1971).}
Domaine d'holomorphie des fonctions v\'erifiant une \'equation aux d\'eriv\'ees partielles.
{\it Comptes Rendus de l'Acad\'emie des Sciences Paris} {\bf 272} (1646--1648).




    \end{thebibliography}
\end{document}